\documentclass[11pt]{amsart}
\usepackage{amsmath,amssymb,amsfonts,url,mathptmx}
\usepackage{color}
\usepackage{appendix}
\allowdisplaybreaks[4]
\numberwithin{equation}{section}

\newtheorem{thm}{Theorem}[section]
\newtheorem{lem}{Lemma}[section]
\newtheorem{rem}{Remark}[section]
\newtheorem{prop}{Proposition}[section]

\newcommand{\hdot}{^\text{\r{}}\hspace{-.33cm}H}
\usepackage[marginal]{footmisc}

\begin{document}

\title[SU(3) Toda system]{Estimates of bubbling solutions of $SU(3)$  Toda systems at critical parameters-Part 2}
\subjclass{}

\author{Juncheng Wei, Lina Wu, Lei Zhang}\footnote{Lei Zhang is partially supported by a Simons Foundation Collaboration Grant}

\address{Department of Mathematics\\
        University of Florida\\
        1400 Stadium Rd\\
        Gainesville FL 32611}
\email{leizhang@ufl.edu}

\date{\today}

\begin{abstract}
In this article we study bubbling solutions of regular $SU(3)$ Toda systems defined on a Riemann surface. There are two major difficulties corresponding to the profile of bubbling solutions: partial blowup phenomenon and bubble accumulation. We prove that when both parameters  tend to critical positions, if there is one fully bubbling blowup point, then under one curvature assumption, all the blowup solutions near a blowup point satisfy a spherical Harnack inequality, which completely rules out the bubble-accumulation phenomenon. This fact is crucial for a number of applications.
\end{abstract}

 \maketitle

\section{Introduction}
In this article we consider the following $SU(3)$ Toda system defined on a compact Riemann surface $(M,g)$:
\begin{align}\label{main-eq-0}
\Delta_g u_1+2\rho_1(\frac{h_1e^{u_1}}{\int_Mh_1e^{u_1}dV_g}-1)-\rho_2(\frac{h_2e^{u_2}}{\int_Mh_2e^{u_2}dV_g}-1)=0\\
\Delta_g u_2-\rho_1(\frac{h_1e^{u_1}}{\int_Mh_1e^{u_1}dV_g}-1)+2\rho_2(\frac{h_2e^{u_2}}{\int_Mh_2e^{u_2}dV_g}-1)=0 \nonumber
\end{align}
where the volume of $M$ is assumed to be $1$ for convenience. $h_1,h_2$ are positive smooth functions, $\rho_1,\rho_2$ are positive constants.

The standard space for $u=(u_1,u_2)$ is that each component is in 
$$\hdot^1(M):=\{v\in H^1(M);\quad \int_M v=0. \}$$ and system (\ref{main-eq-0}) is the Euler Lagrange equation of the following variational form:
$$J_{\rho}(u)=\frac 12\int_M\sum_{i,j=1}^2a^{ij}\nabla_gu_i\nabla_gu_jdV_g-\sum_{i=1}^2\rho_i\log\int_M h_ie^{u_i}dV_g $$
where $a^{11}=a^{22}=\frac 23$, $a^{12}=a^{21}=\frac 13$. It is proved by Lin-Wei-Yang-Zhang in \cite{lin-wei-yang-zhang-apde} that a priori estimate holds if $(\rho_1,\rho_2)\in (4\pi (m-1),4\pi m)\times (4\pi (n-1),4\pi n)$ for $m,n\in \mathbb N$ (the set of natural numbers). Thus a degree-counting program can be defined by this a priori estimate. Since this degree counting program links the local blowup analysis with the topology of $M$, it reveals important interplay of different fields of mathematics. 

Even though the a priori estimate requires $\rho=(\rho_1,\rho_2)$ to be away from multiples of $4\pi$, it is also important to study blowup solutions when both parameters tend to multiples of $4\pi$. For comparison if the Toda system is reduced to one equation, which is the Liouville equation, the critical parameter for Liouville equation is when the parameter $\rho$ is a multiple of $8\pi$. The case $\rho=8\pi$ is studied by Lin-Wang \cite{lin-wang-ann} that reveals significant ties among elliptic functions, blowup solutions and the number of critical points of the Green's function on a torus. The study of Toda system with critical parameters should also lead to possibly surprising connections between fields. In fact (\ref{main-eq-0}) is the simplest Toda system that connects many fields of mathematics and physics. Since the literature is too vast to be listed, we just mention the following works closely related to the topic in this article \cite{ao-wang,bart-taran-mass,bart-nonsimple,chen-lin-last-cpam,chen-lin-wang,april-pistoia,del-pino-1,del-pino-2,jevnikar,karmakar,kuo-lin,lee-4-1,lee-4-2,lin-wei-yang-zhang-apde,lin-wei-zhang-adv,lin-wei-zhao,lin-yang,lin-yang-zhong}. As a matter of fact, the research of Toda system has been so dynamical that many breakthroughs come from inspirations of different fields and exciting ties among seemingly far-reached fields have been greatly strengthened by new discoveries. 

A sequence of solutions $u^k=(u_1^k,u_2^k)$ is called  bubbling solutions, or blowup solutions, if $\max_{M}\{u_1^k(p_k),u_2^k(p_k)\}\to \infty$ for $p_k\to p\in M$. If $u^k$ converges to a global $SU(3)$ Toda system after scaling according to the maximum of both components, $u^k$ is called a \emph{fully bubbling } sequence. Otherwise $u^k$ is called to have a \emph{partial blowup} phenomenon. In other words, fully bubbling solutions converge to an entire Toda system, while partially blown-up solutions tend to only one equation after scaling. In blowup analysis there are two major difficulties, the first one is the partial blowup phenomenon, the second one is the bubble-accumulation situation: disjoint bubbling disks may tend to the same point, which makes the study of the profile of bubbling solutions particularly challenging. 

The main purpose of this article is to remove the bubble accumulation phenomenon and thus greatly simplifies the profile of bubbling solutions.  The following is a curvature assumption we make on the coefficient function $h_i^k$: For each $x\in M$, 
\begin{equation}\label{curvature-a}
\Delta_g (\log h_i)(x)-2K(x)\not \in 4\pi \mathbb Z,\quad i=1,2.
\end{equation}
where $\mathbb Z$ stands for the set of integers, $K(x)$ is the Gauss curvature at $x$. Then our main theorem is
\begin{thm}\label{main-thm} Let $u^k$ be a sequence of blowup solutions to (\ref{main-eq-0}) with parameters $\rho^k\to (4\pi m,4\pi n)$ for some positive integers $m$ and $n$. If there exists a fully bubbling blowup point and (\ref{curvature-a}) holds, 
around each blowup point $p_l$ we have 
\begin{equation}\label{main-sphe-har}
u_i^k(x)+2\log |x-p_l^k|\le C, \quad i=1,2, 
\end{equation}
where $p_l^k\to p_l$ is a local maximum of one component of $u^k$.
\end{thm}
The inequality (\ref{main-sphe-har}) implies that around $p_l^k$, the spherical Harnack inequality holds for both components, thus Theorem \ref{main-thm} completely rules out the bubble-accumulation phenomenon and provides a key estimate for construction of bubbling solutions as well the degree-counting program mentioned before. It is also important to point out that Theorem \ref{main-thm} does not rule out partial blowup phenomenon. The construction of a $(4,2)$ type blowup by Musso-Pistoia-Wei \cite{musso-pistoia-wei}, which contains a partial blowup situation, does satisfy the spherical Harnack inequality around a blowup point. It is an interesting open question whether (\ref{curvature-a}) also leads to the elimination of partial blowup phenomenon. Here we also point out that having a simple blowup situation leads to key information about curvature as well. For example, in \cite{lin-wei-zhao}, \cite{lin-wei-zhang-adv} and \cite{zhang-imrn}, only fully bubbling solutions are assumed for Toda systems, the authors proved surprising vanishing estimates on coefficient functions, which serve as guidelines of bubble construction, degree counting and other applications. 

\begin{rem}
The assumption $\rho^k=(\rho_1^k,\rho_2^k)\to (4m\pi,4n\pi)$ is more general than the assumption in  \cite{wu-zhang-part-1} where $\rho_k\to (4\pi, 4m\pi)$ or $\rho^k\to (4m\pi,4\pi)$ was assumed.
\end{rem}
Since the proof of the main theorem is quite involved, we describe an outline of the proof at the end of the introduction. The statement of Theorem \ref{main-thm} is equivalent to saying that the spherical Harnack inequality holds around a local maximum of a component. As is well known, if the spherical Harnack inequality does not hold, one would obtain bubbling disks mutually disjoint from one another. Inside each bubbling disk, the profile of the blowup solutions is either that of a global solution of the $SU(3)$ Toda system, or that of a global solution of a single Liouville equation. To eliminate the bubble accumulation we employ the ideas in a series of paper of Wei-Zhang \cite{wei-zhang-adv}, \cite{wei-zhang-plms} and \cite{wei-zhang-jems}. One major difficulty is to rule out the case that one bubbling disk of $u_1^k$ is placed between two bubbling disks of $u_2^k$. This part is placed in section three in three stages: In stage one: After removing $u_1^k$ from the equation of $u_2^k$, the equation of $u_2^k$ is similar to a singular Liouville equation with a quantized singular source. $u_2^k$ has two local maximum points and there is a Pohozaev identity corresponding to each local maximum. After working on these two Pohozaev identities we shall see that the first derivatives of the coefficient functions are involved with the location of blowup points. Next we compare the two pohozaev identities with those of a global solution. The comparison of Pohozaev identities will lead to an initial vanishing rate of certain coefficient functions. Then in stage two we improve the pointwise estimate as well as the vanishing rate of the first derivative of the coefficient functions. The key point of the proof of stage two is that the blowup solutions look like \emph{different global solutions } around each local maximum, which is the main reason of the contradiction. The in stage three, we employ the argument of stage two and the vanishing rate of the first derivatives of the coefficient functions to obtain a better global pointwise estimate. In this case the precise pointwise estimate of the single Liouville equation plays a crucial role. These three stages rule out the case $(2,4)$ as well as $(4,2)$ under the curvature condition stated in Theorem \ref{main-thm}. The other cases of blowup solutions violating the spherical harnack inequalities are relatively easier and they are handled in section four. There are quite a few key estimates in the proof and they are explained in details as the proof is carried out. 

\section{A locally defined system}
In this section we write the $SU(3)$ Toda system as a set of locally defined equations near a fixed point. 
For any fixed point we use isothermal coordinates: Let $g_0$ be the flat metric and $\phi$ satisfy $g_0=e^{-\phi}g$. Then $\Delta_g=e^{-\phi}\Delta$ and $\phi$ satisfies
\begin{equation}\label{local-phi}
\phi(0)=|\nabla \phi(0)|=0,\quad \Delta \phi(0)=-2K(p). 
\end{equation}
Let $\hat u_i=u_i-\log \int_M h_ie^{u_i},\quad i=1,2$. Then we have $\int_M h_ie^{\hat u_i}=1$. 
Then for $\hat u=(\hat u_1,\hat u_2)$ we have
\begin{align*}
\Delta \hat u_1+2\rho_1 h_1e^{\hat u_1+\phi}-\rho_2h_2e^{\hat u_2+\phi}=(2\rho_1-\rho_2)e^{\phi},\\
\Delta \hat u_2-\rho_1 h_1 e^{\hat u_1+\phi}+2\rho_2h_2e^{\hat u_2+\phi}=(2\rho_2-\rho_1)e^{\phi}
\end{align*}
Let $B_{\tau}$ be a ball of radius $\tau$ centered at the origin. For $\tau>0$ small we set $f_1$ and $f_2$ as
\begin{align}\label{local-f}
\Delta f_1=(2\rho_1-\rho_2)e^{\phi},\quad f_1(0)=0,\quad f_1=\mbox{constant on } \partial B_{\tau} \\
\Delta f_2=(2\rho_2-\rho_1)e^{\phi},\quad f_2(0)=0,\quad f_2=\mbox{constant on }\partial B_{\tau}. \nonumber
\end{align}
Then by setting 
$$\mathfrak{u}_i=\hat u_i-f_i, \quad i=1,2 $$
and 
$$\mathfrak{h}_i=\rho_ih_ie^{\phi+f_i},\quad i=1,2 $$
in a neighborhood of $0$, say $B_{\tau}$ we have
\begin{align}\label{main-local}
   & \Delta \mathfrak{u}_1+2\mathfrak{h}_1e^{\mathfrak{u}_1}-\mathfrak{h}_2e^{\mathfrak{h}_2}=0,\\
   & \Delta \mathfrak{u}_2-\mathfrak{h}_1e^{\mathfrak{u}_1}+2\mathfrak{h}_2e^{\mathfrak{h}_2}=0,\quad \mbox{in}\quad B_{\tau}. \nonumber
\end{align}

\section{Blowup analysis for a local system}
In this long section we derive a vanishing estimate for a locally defined Toda system.
Let $B_1$ be the unit ball in $\mathbb R^2$ and we consider the following locally defined system:
\begin{align}\label{main-eq}
&\Delta \mathfrak{u}_1^k+ 2h_1^ke^{\mathfrak{u}_1^k}-h_2^ke^{\mathfrak{u}_2^k}=0,\nonumber\\
&\Delta \mathfrak{u}_2^k-h_1^ke^{\mathfrak{u}_1^k}+2h_2^ke^{\mathfrak{u}_2^k}=0,\quad \mbox{in}\quad B_1,
\end{align}
and postulate usual assumptions:
\begin{equation}\label{bound-1}
\max_{x\in K}\{\mathfrak{u}_1^k(x),\mathfrak{u}_2^k(x)\}\le C(K),\quad K\subset\subset B_1\setminus\{0\}.
\end{equation}
Boundary oscillation finite-ness:
\begin{equation}\label{boundary-f}
\max |\mathfrak{u}_i^k(x)-\mathfrak{u}_i^k(y)|\le C, \quad i=1,2,\quad \forall x,y\in \partial B_1.
\end{equation}
and there exists a $C>0$ such that
\begin{equation}\label{assum-h}
\frac{1}{C}\le h_i^k(x)\le C, \quad \|h_i^k\|_{C^3(B_1)}\le C. 
\end{equation}
and there is a uniform bound on the total integration:
\begin{equation}\label{finite-energy}
\int_{B_1}h_i^ke^{\mathfrak{u}_i^k}\le C. 
\end{equation}

Let $\phi_i^k$ be the harmonic function defined by the oscillation of $\mathfrak{u}_i$ on $\partial B_1$:
$$\Delta \phi_i^k(x)=0,\quad \mbox{ in } B_1, \qquad \phi_i^k(x)=\mathfrak{u}_i^k(x)-\frac{1}{2\pi}\int_{\partial B_1}\mathfrak{u}_i^kdS,\quad x\in \partial B_1. $$

The following quantity is used to describe the concentration of the integrals of $e^{\mathfrak{u}_i^k}$:
 Let 
$$\sigma_i=\frac{1}{2\pi}\lim_{\delta\to 0}\lim_{k\to \infty}\int_{B_{\delta}}h_i^ke^{\mathfrak{u}_i^k}. $$ 
Note that in this article we don't distinguish sequence and a sub-sequence.
Then it is well known (see \cite{wu-zhang-part-1}, for example) 
that
$$(\sigma_1,\sigma_2)=(2,0),(0,2)(2,4),(4,2),(4,4). $$
The main result of this section is:

\begin{thm}\label{local-laplace}
Let $\mathfrak{u}^k=(\mathfrak{u}_1^k,\mathfrak{u}_2^k)$ be solutions of (\ref{main-eq}) such that (\ref{bound-1}),(\ref{boundary-f}), (\ref{assum-h}) and (\ref{finite-energy}) hold for $\mathfrak{h}_i$ and $\mathfrak{u}^k$. Suppose $\mathfrak{u}^k$ is a $(2,4)$ type blowup sequence that violates the spherical Harnack inequality around the local maximum of $\mathfrak{u}_1^k$, and 
\begin{equation}\label{mu-1-large}
\max_{B_1}\mathfrak{u}_1^k>\max_{B_1}\mathfrak{u}_2^k+5\log \delta_k^{-1}
\end{equation}
where $\delta_k$ is the distance between one local maximum of $\mathfrak{u}_2^k$ and one local maximum of $u_1^k$. Then we have
$$\nabla \log h_2^k(0)+\nabla \phi_2^k(0)=o(1), $$
and
$$\Delta \log h_2^k(0)=o(1) $$
\end{thm}
\begin{rem} There are two formations of the blowup type $(2,4)$. The one that violates the spherical Harnack inquality has one bubbling disk of $\mathfrak{u}_1^k$ in the middle of two bubbling disks of $\mathfrak{u}_2^k$ on both sides. The second formation consists of a bubbling disk of $\mathfrak{u}_1^k$ on top of $\mathfrak{u}_2^k$. In the second case, if $\mathfrak{u}_2^k$ is scaled according to its maximum, the limit function is 
$$\Delta v+2 e^{v}=4\pi \delta_0,\quad \mbox{in}\quad \mathbb R^2, \quad \int_{\mathbb R^2}e^{v}<\infty. $$
What Theorem \ref{local-laplace} covers is the first case. 
\end{rem}

\noindent{\bf Proof of Theorem \ref{local-laplace}:} 

Let $p_0^k$ be the location of the local maximum of $\mathfrak{u}_1^k$, $p_1^k$ and $p_2^k$ be the locations of local maximums if $\mathfrak{u}_2^k$. In the $(2,4)$ formation, it is known that \cite{wu-zhang-part-1}
$$\lim_{k\to \infty}\frac{|p_1^k-p_0^k|}{|p_1^k-p_2^k|}=\frac 12, \quad 
\lim_{k\to \infty}\frac{p_1^k-p_0^k}{|p_1^k-p_0^k|}=\lim_{k\to \infty}\frac{p_1^k-p_2^k}{|p_1^k-p_2^k|}. $$
Suppose $\delta_k=|p_1^k-p_0^k|$ and we write $p_1^k=p_0^k+\delta_ke^{i\theta_k}$. Then we set 
$$\tilde v_i^k(y)=\mathfrak{u}_i^k(p_0^k+\delta_ke^{i\theta_k}y)+2\log \delta_k, \quad i=1,2. $$

Now the system can be written as
\begin{align}\label{local-eq-1}
&\Delta \tilde v_1^k(y)+2h_1^k(p_0^k+\delta_ky)e^{\tilde v_1^k}-h_2^k(p_0^k+\delta_ky)e^{\tilde v_2^k}=0, \\
&\Delta \tilde v_2^k(y)-h_1^k(p_0^k+\delta_ky)e^{\tilde v_1^k}+2h_2^k(p_0^k+\delta_ky)e^{\tilde v_2^k}=0,\quad 
\mbox{in}\quad \Omega_k\nonumber
\end{align}
where $\Omega_k=B(0,\tau\delta_k^{-1})$ for some $\tau>0$ and $G_k$ be the Green's function on $\Omega_k$ with respect to the Dirichlet boundary condition: 
$$G_k(y,\eta)=-\frac 1{2\pi}\log |y-\eta |+H_k(y,\eta)$$
where 
\begin{equation}\label{Hk-exp}
H_k(y,\eta)=\frac 1{2\pi}\log (\frac{|\eta |}{L_k}|\frac{L_k^2\eta}{|\eta |^2}-y|),\quad L_k=\tau \delta_k^{-1}. 
\end{equation}
Then we remove $\tilde v_1^k$ from the equation of $\tilde v_2^k$.
Let 
$$f_1^k=-\int_{\Omega_k}G_k(y,\eta)h_1^ke^{\tilde v_1^k}d\eta+\int_{\Omega_k}G_k(0,\eta)h_1^ke^{\tilde v_1^k}d\eta$$
we have 
$$\Delta f_1^k=h_1^k(p_0^k+\delta_ky)e^{\tilde v_1^k},\quad f_1^k(0)=0$$
and let 
$$v_2^k=\tilde v_2^k-f_1^k,$$
and
$$\mathfrak{h}_i^k(\delta_ky):=h_i^k(p_0^k+\delta_k y), \quad i=1,2. $$
Since we will use complex numbers to represent points on $\mathbb R$, we use $e_1$ to denote $(1,0)$ and $e^{i\pi}$ to denote $(-1,0)$. For $y$ around $e_1$ or $e^{i\pi}$, 
since the singularity of $\tilde v_1^k$ is at the origin, standard estimate can be used to obtain
$$f_1^k(y)=2\log |y|+O(\mu_1^ke^{-\mu_1^k})+O(\delta_k^2),\quad y\in B(e_1,\tau)\cup B(e^{i\pi},\tau) $$
for some $\tau>0$ small, where we use $\mu_1^k$ and $\mu_2^k$ to denote the heights of the bubbles:
$$\mu_1^k=\max v_1^k, \quad \mu_2^k=v_2^k(e_1). $$
In this section we use $E_1$ to denote
$$E_1=O(\mu_1^ke^{-\mu_1^k})+O(\delta_k^2). $$
Then the equation of $v_2^k$ now becomes
\begin{equation}\label{v2-eq-weak}
\Delta v_2^k+|y|^{2+E_1}\mathfrak{h}_2^k(\delta_ky)e^{v_2^k}=0 
\end{equation}

It is important to observe that after scaling, the local maximum of $v_1^k$ is the origin and one local maximum of $v_2^k$ is $e_1$, the other local maximum of $v_2^k$ is $e^{i\pi}+o(1)$. We use $Q_0^k$ to denote $e_1$ and $Q_1^k$ to denote the other local maximum of $v_2^k$. (add a footnote)

By the pointwise analysis of single Liouville equation, $v_2^k$ is around $-\mu_2^k$ away from bubbling disks in $B_{10}$, $v_1^k(y)=-\mu_1^k+O(1)$ in the same area. 

Now we provide a rough estimate for $v_i^k$ outside the bubbling area:

For $r>2$, let $\bar v_2^k(r)$ be the spherical average of $v_2^k$ on $\partial B_r$, then we have
$$\frac{d}{dr}\bar v_2^k(r)=\frac{d}{dr}\bigg (\frac{1}{2\pi r}\int_{B_r}\Delta v_2^k\bigg )=-\frac{12\pi+o(1)}{2\pi r}. $$
Because of the fast decay of $\bar v_2^k(r)$ it is easy to use the Green's representation of $v_2^k$ to obtain the following stronger estimates:
\begin{align}\label{vk-crude}
v_2^k(y)&=- \mu_2^k-(6+o(1))\log |y|+O(1),\quad 2<|y|<\tau \delta_k^{-1}.\nonumber\\
v_1^k(y)&=-\mu_1^k+o(1)\log |y|+O(1),\quad 2<|y|<\tau \delta_k^{-1}. 
\end{align}

\subsection{An initial vanishing estimate of the first derivatives}

Now we consider $v_2^k$ around $Q_l^k$. Using the results in \cite{chenlin1,zhangcmp,gluck} we have, for $v_2^k$ in $B(Q_l^k,\epsilon )$, the following gradient estimate:

 \begin{equation}\label{gra-each-p}
\delta_k \nabla  (\log \mathfrak{h}_2^k)(\delta_k \tilde Q_l^k)+2\frac{\tilde Q_l^k}{|\tilde Q_l^k|^2}+\nabla \phi_l^k(\tilde Q_l^k)=O(\mu_2^ke^{-\mu_2^k}),\quad l=0, 1
\end{equation}
where $\phi_l^k$ is the harmonic function that eliminates the oscillation of $v_2^k$ on $\partial B(Q_l^k,\epsilon)$ and $\tilde Q_l^k$ is the maximum of $v_2^k-\phi_l^k$ that satisfies
\begin{equation}\label{close-2}
\tilde Q_l^k-Q_l^k=O( e^{-\mu_2^k}).
\end{equation}

Using (\ref{close-2}) in (\ref{gra-each-p}) we have
\begin{equation}\label{pi-each-p}
\delta_k \nabla  (\log \mathfrak{h}_2^k)(\delta_k  Q_l^k)+2\frac{Q_l^k}{| Q_l^k|^2}+\nabla \phi_l^k( Q_l^k)=O(\mu_2^ke^{-\mu_2^k}).
\end{equation}
For the discussion in this section we use the following version of (\ref{pi-each-p}):
\begin{equation}\label{location-delta}
\delta_k\nabla  (\log \mathfrak{h}_2^k)(0)+2\frac{Q_l^k}{|Q_l^k|^2}+\nabla \phi_l^k(Q_l^k)=O(\delta_k^2)+O(\mu_2^ke^{-\mu_2^k})
\end{equation}
and the first estimate of $\nabla \phi_l^k(Q_l^k)$ is
\begin{lem}\label{phi-k-e} For $l=0,1$,
\begin{equation}\label{gra-each-2}
\nabla \phi_l^k(Q_l^k)
=-4\frac{Q_l^k-Q_m^k}{|Q_l^k-Q_m^k|^2}
+O(\mu_2^ke^{-\mu_2^k}) \quad l\neq m.
\end{equation}
\end{lem}

\noindent{\bf Proof of Lemma \ref{phi-k-e}:}

From the expression of $v_2^k$ on $\Omega_k=B(0,\tau \delta_k^{-1})$ we have, for $y$ away from bubbling disks,
\begin{align}\label{pi-e-1}
v_2^k(y)&=v_2^k|_{\partial \Omega_k}+\int_{\Omega_k}G_k(y,\eta)|\eta |^{2}\mathfrak{h}_2^k(\delta_k\eta)e^{v_2^k(\eta)}d\eta\\
&=v_2^k|_{\partial \Omega_k}+\sum_{l=0}^{1}G_k(y,Q_l^k)\int_{B(Q_l,\epsilon)}|\eta |^{2}\mathfrak{h}_2^k(\delta_k\eta)e^{v_2^k}d\eta \nonumber\\
&+\sum_l\int_{B(Q_l,\epsilon)}(G_k(y,\eta)-G_k(y,Q_l^k))|\eta |^{2}\mathfrak{h}_2^k(\delta_k\eta)e^{v_2^k}d\eta+O(\mu_2^ke^{-\mu_2^k}). \nonumber
\end{align}
Before we evaluate each term, we use a sample computation which will be used repeatedly: Suppose $f$ is a smooth function defined on $B(Q_0^k,\epsilon)$, then
we evaluate
\begin{equation}\label{sample-c}
\int_{B(Q_0^k,\epsilon)}f(\eta)|\eta |^{2}\mathfrak{h}_2^k(\delta_k \eta)e^{v_2^k(\eta)}d\eta.
\end{equation}
Let $\tilde Q_0^k$ be the maximum of $v_2^k-\phi_0^k$, and set $\hat h_k(y)=|y|^{2}\mathfrak{h}_2^k(\delta_ky)$ then it is known \cite{zhangcmp,gluck} that
\begin{equation}\label{close-Q}
\tilde Q_0^k-Q_0^k=O(e^{-\mu_2^k}).
\end{equation}
Moreover, it is derived that
\begin{align}\label{gluck-t}
v_2^k(y)=\phi_0^k(y)+\log \frac{e^{\mu_2^k}}{(1+e^{\mu_2^k}\frac{\hat h_k(\tilde Q_0^k)}{8}|y-\tilde Q_0^k|^2)^2}\\
+c_ke^{-\mu_2^k}(\log (2+e^{\mu_2^k/2}|y-\tilde Q_0^k|))^2+O(\mu_2^ke^{-\mu_2^k}) \nonumber
\end{align}
for some $c_k\in \mathbb R$. We use $U_k$ to denote the leading global solution.
Using the expansion above and symmetry we have
\begin{align}\label{sampled-cm}
&\int_{B(\tilde Q_0^k,\epsilon)}f(\eta)|\eta |^{2}\mathfrak{h}_2^k(\delta_k \eta)e^{v_2^k(\eta)}d\eta\\
=&\int_{B(\tilde Q_0^k,\epsilon)}f(\eta)\hat h_k(y)e^{v_2^k(y)}dy\nonumber \\
=&\int(f(\tilde Q_0^k)+\nabla f(\tilde Q_0^k)\tilde \eta+O(|\tilde \eta|^2)(\hat h_k(\tilde Q_0^k)+\nabla \hat h_k(\tilde Q_0^k)(\tilde \eta)
+O(|\tilde \eta|^2)\nonumber \\
&\cdot e^{U_k}(1+\phi_0^k+c_ke^{-\mu_2^k}(\log (2+e^{\mu_2^k/2}|\tilde \eta|))^2+O(\mu_2^ke^{-\mu_2^k})+O(|\tilde \eta|^2))d\eta \nonumber
\end{align}
where $\tilde \eta=\eta-\tilde Q_0^k$.  In the evaluation we use symmetry, for example,
$$\int_{B(\tilde Q_0^k,\epsilon)}e^{U_k}\tilde \eta d\eta=0. $$ Also
using $\phi_0^k(Q_0^k)=0$ and $Q_0^k-\tilde Q_0^k=O(e^{-\mu_2^k})$ we have
$$\int_{B(\tilde Q_0^k,\epsilon)} e^{U_k}\phi_0^k=O(e^{-\mu_k}). $$
Carrying out computations in (\ref{sampled-cm}) we arrive at
$$\int_{B(Q_0^k,\epsilon)}f(\eta)|\eta |^{2}\mathfrak{h}_2^k(\delta_k \eta)e^{v_2^k(\eta)}d\eta=8\pi f(\tilde Q_0^k)+O(\mu_2^k e^{-\mu_2^k}). $$
Since (\ref{close-Q}) holds we further have
\begin{equation}\label{sample-2}
\int_{B(Q_0^k,\epsilon)}f(\eta)|\eta |^{2}\mathfrak{h}_2^k(\delta_k \eta)e^{v_2^k(\eta)}d\eta=8\pi f( Q_0^k)+O(\mu_2^k e^{-\mu_2^k}).
\end{equation}

Using the method of (\ref{sample-2}) in the evaluation of each term in (\ref{pi-e-1}) we have,
$$
v_2^k(y)
=v_2^k|_{\partial \Omega_k}-4\sum_{l=0}^{1}\log |y- Q_l^k|
+8\pi\sum_{l=0}^{1}H_k(y,Q_l^k)+O(\mu_2^ke^{-\mu_2^k}). $$

The harmonic function that kills the oscillation of $v_2^k$ around $Q_m^k$ is
\begin{align*}
\phi_m^k&=-4(\log |y-Q_l^k|-\log |Q_m^k-Q_l^k|)\\
&+8\pi\sum_{s=0}^{1}(H_k(y,Q_s^k)-H_k(Q_m^k,Q_s^k))+E,\quad l\neq m.
\end{align*}

The corresponding estimate for $\nabla \phi_m^k$ is
\begin{equation}\label{Phk-exp}
\nabla \phi_m^k(Q_m^k)=-4\frac{Q_m^k-Q_l^k}{|Q_m^k-Q_l^k|^2}
+8\pi\sum_{s=0}^{1}\nabla_1 H_k(Q_m^k,Q_s^k)+O(\mu_2^ke^{-\mu_2^k}),\,l\neq m,
\end{equation}
where $\nabla_1$ stands for the differentiation with respect to the first component. From the expression of $H_k$ in (\ref{Hk-exp}), we have
\begin{align}\label{grad-H}
\nabla_1H_k(Q_m^k,Q_l^k)&=\frac 1{2\pi}\frac{Q_m^k-\tau^2\delta_k^{-2}Q_l^k/|Q_l^k|^2}{|Q_m^k-\tau^2\delta_k^{-2}Q_l^k/|Q_l^k|^2|^2}\\
&=\frac 1{2\pi}\tau^{-2}\delta_k^2\frac{\tau^{-2}\delta_k^2Q_m^k-Q_l^k/|Q_l^k|^2}{|Q_l^k/|Q_l^k|^2-\tau^{-2}\delta_k^2Q_m^k|^2}\nonumber \\
&=-\frac{1}{2\pi}\tau^{-2}\delta_k^2e^{i\pi l}+O(\sigma_k\delta_k^2). \nonumber
\end{align}
where $\sigma_k=\max_l|Q_l^k-e^{\pi i l}|$. Later we shall obtain more specific estimate of $\sigma_k$.

Thus using (\ref{grad-H}) in (\ref{Phk-exp}) we have 
\begin{align}\label{phimk-g}
&\nabla \phi_m^k(Q_m^k) \\
=&-4\frac{Q_m^k-Q_l^k}{|Q_m^k-Q_l^k|^2}
-4\tau^{-2}\delta_k^2\sum_{l=0}^{1}e^{\pi il}+O(\sigma_k\delta_k^2)+O(\mu_ke^{-\mu_k}) \nonumber \\
=&-4\frac{Q_m^k-Q_l^k}{|Q_m^k-Q_l^k|^2}
+O(\sigma_k\delta_k^2)+O(\mu_2^ke^{-\mu_2^k})\quad l\neq m.\nonumber
\end{align}
Since we don't have the estimate of $\sigma_k$ we use the following weak version for now:
$$
\nabla \phi_m^k(Q_m^k)=-4\frac{Q_m^k-Q_l^k}{|Q_m^k-Q_l^k|^2}+E, \quad l\neq m.
$$
where
$$E=O(\delta_k^2)+O(\mu_2^ke^{-\mu_2^k}). $$
Lemma \ref{phi-k-e} is established. $\Box$

\subsection{Location of blowup points}

In this subsection we establish a first description of the locations of $Q_1^k$. The result of this section will be used later to obtain vanishing estimates of the coefficient function $h_2^k$.

The Pohozaev identity around $Q_1^k$ now reads
\begin{equation}\label{tem-pi}
-4\frac{Q_l^k-Q_j^k}{|Q_l^k-Q_j^k|^2}+2\frac{Q_l^k}{|Q_l^k|^2}=-\nabla (\log \mathfrak{h}_2^k)(0)\delta_k +E,\quad j\neq l. 
\end{equation}
Using $L$ to denote $\nabla (\log  \mathfrak{h}_2^k)(0)$:
$$L=\nabla (\log  \mathfrak{h}_2^k)(0),$$
we write (\ref{tem-pi}) as 
$$\frac{1}{\bar Q_l^k}=2\frac{1}{\bar Q_l^k-\bar Q_j^k}-\frac{L}2\delta_k+E, j\neq l $$
where $\bar Q_l^k$ is the conjugate of $Q_l^k$.

$Q_0=e_1$, $Q_1=-e_1-m_1^k$ where every number is in the complex plane. Now the two equations from the Pohozaev identity are
\begin{eqnarray*}
1-\frac{2}{2+\bar m_1^k}=-\frac{L}{2}\delta_k+E\\
-\frac{1}{1+\bar m_1^k}+\frac{2}{2+\bar m_1^k}=-\frac{L}2\delta_k+E. 
\end{eqnarray*}

From these two equations we obtain
\begin{equation}\label{m-2}
\bar m_1^k=-L\delta_k+E.
\end{equation}

With this fact we can further write $\nabla_1 H_k(Q_m^k,Q_l^k)$  in (\ref{grad-H}) as
\begin{equation}\label{grad-H-1}
\nabla_1H_k(Q_m^k,Q_l^k)
=-\frac{1}{2\pi}\tau^{-2}\delta_k^2e^{\pi i l}+O(\delta_k^3)+O(\mu_2^ke^{-\mu_2^k}),
\end{equation}
and $\nabla \phi_l^k(Q_l^k)$ in (\ref{gra-each-2}) and (\ref{phimk-g}) as
\begin{equation}\label{gra-each-3}
\nabla \phi_l^k(Q_l^k)
=-4\frac{Q_l^k-Q_m^k}{|Q_l^k-Q_m^k|^2}
+O(\delta_k^3)+O(\mu_2^ke^{-\mu_2^k}),\quad m\neq l.
\end{equation}

\subsection{Initial Vanishing estimates of $\nabla \mathfrak{h}_2^k(0)$}

 Before we get into the details we would like to mention that the most important observation is that the difference between the Pohozaev identities of $v_2^k$ and that of a global solution can be evaluated in detail. Using the result in the previous section we shall see that the coefficient of  $\nabla\mathfrak{h}_2^k(0)$ is not zero.

First we recall the equation for $v_2^k$:
$$\Delta v_2^k+\mathfrak{h}_2^k(\delta_ky)|y|^{2}e^{v_2^k}=E_1,\quad |y|<\tau \delta_k^{-1} $$
with $v_2^k=$ constant on $\partial B(0,\tau \delta_k^{-1})$. Moreover $v_2^k(e_1)=\mu_2^k$. Now we set
$V_k$ to be the solution to
$$\Delta V_k+\mathfrak{h}_2^k(\delta_ke_1)|y|^{2}e^{V_k}=0,\quad \mbox{in}\quad \mathbb R^2, \quad \int_{\mathbb R^2}|y|^{2}e^{V_k}<\infty $$
such that $V_k$ has its local maximums at $e_1$ and $-e_1$, and $V_k(e_1)=\mu_2^k$.
The expression of $V_k$ is
\begin{equation}\label{Vk-exp}
V_k(y)=\log \frac{e^{\mu_2^k}}{(1+\frac{e^{\mu_2^k}\mathfrak{h}_2^k(\delta_ke_1)}{32}|y^{2}-e_1|^2)^2}.
\end{equation}
In particular for $|y|\sim \delta_k^{-1}$,
$$V_k(y)=-\mu_2^k-8\log \delta_k^{-1}+C+O(\delta_k^{2})+O(e^{-\mu_2^k}). $$

Let $w_k=v_2^k-V_k$ and $\Omega_k=B(0,\tau \delta_k^{-1})$,  we shall derive a precise, point-wise estimate of $w_k$ in $B_3\setminus B(Q_1,\tau)$ where $\tau>0$ is a small number independent of $k$. Here we note that among $2$ local maximum points, we already have $e_1$ as a common local maximum point for both $v_2^k$ and $V_k$ and we shall prove that $w_k$ is very small in $B_3$ if we exclude $B(Q_1,\tau)$. Before we carry out more specific computations we emphasize the importance of
\begin{equation}\label{control-e}
w_k(e_1)=|\nabla w_k(e_1)|=0.
\end{equation}
First we consider the Green's representation of $v_2^k$ on $\Omega_k$
$$v_2^k(y)=\int_{\Omega_k}G_k(y,\eta)|\eta|^{2}\mathfrak{h}_2^k(\delta_k\eta)e^{v_2^k(\eta)}d\eta+v_2^k|_{\partial \Omega_k}, $$
where in $B(0,L_k)$ ($L_k=\tau \delta_k^{-1}$).
Similarly for $V_k$ we obtain from the asymptotic expansion of $V_k$ to have
$$V_k(y)=\int_{\Omega_k}G_k(y,\eta)|\eta|^{2}\mathfrak{h}_2^k(\delta_ke_1)e^{V_k(\eta)}d\eta+\bar V_k|_{\partial \Omega_k}+O(\delta_k^{2})+O(e^{-\mu_2^k}), $$
where $\bar V_k|_{\partial \Omega_k}$ is the average of $V_k$ on $\partial \Omega_k$ and we have used the fact that the oscillation of $V_k$ on $\partial \Omega_k$ is $O(\delta_k^{2})+O(e^{-\mu_2^k})$.
The combination of these two equations gives, for $y\in B_3$,
$$w_k(y)=\int_{\Omega_k}G_k(y,\eta)|\eta|^{2}(\mathfrak{h}_2^k(\delta_k\eta)e^{v_k}-\mathfrak{h}_2^k(\delta_ke_1)e^{V_k})d\eta+c_k+O(\delta_k^2)+O(e^{-\mu_2^k}). $$
Using $w_k(e_1)=0$, we clearly have
\begin{align}\label{wk-pre-1}
w_k(y)=\int_{\Omega_k}(G_k(y,\eta)-G_k(e_1,\eta))|\eta |^{2}(\mathfrak{h}_2^k(\delta_k\eta)e^{v_k}-\mathfrak{h}_2^k(\delta_ke_1)e^{V_k})d\eta\nonumber\\
+O(\delta_k^2)+O(\mu_2^ke^{-\mu_2^k}).
\end{align}
If we concentrate on $y\in B_3\setminus \cup_{l=0}^{1}B(Q_l^k,\tau)$ for $\tau>0$ small, we first claim that the harmonic part of $G_k$ only contributes $O(\delta_k^2)$ in the estimate:
\begin{equation}\label{har-minor}
\int_{\Omega_k}(H_k(y,\eta)-H_k(e_1,\eta))|\eta |^{2}(\mathfrak{h}_2^k(\delta_k\eta)e^{v_k}-\mathfrak{h}_2^k(\delta_ke_1)e^{V_k})d\eta=O(\delta_k^2)
\end{equation}
where $G_k(y,\eta)=-\frac 1{2\pi}\log |y-\eta |+H_k(y,\eta)$. The proof of (\ref{har-minor}) is elementary and does not require delicate estimates: From
$$H_k(y,\eta)=\frac 1{2\pi}\log L_k+\frac{1}{2\pi}\log |\frac{y}{|y|}-\frac{\eta |y|}{L_k^2}|,\quad L=\tau\delta_k^{-1},$$
we see that for $y\in B_3$,
$$|H_k(y,\eta)-H_k(e_1,\eta)|\le CL_k^{-2}|\eta |. $$
Thus if $|\eta |$ is $O(1)$, the desired estimate is trivial. For $|\eta |$ large, we need to use the fast decay of $e^{v_k}$ and $e^{V_k}$ to prove that the corresponding integral is at most $O(L_k^{-2})$. Thus (\ref{har-minor}) holds.

Combining (\ref{wk-pre-1}) and (\ref{har-minor}), we have
\begin{align}\label{crucial-1}
w_k(y)=\int_{\Omega_k}\frac 1{2\pi}(\log |e_1-\eta|-\log |y-\eta|)|\eta |^{2}(\mathfrak{h}_2^k(\delta_k\eta)e^{v_2^k}-\mathfrak{h}_2^k(\delta_ke_1)e^{V_k})d\eta \nonumber\\
+O(\delta_k^2)+O(e^{-\mu_2^k}).
\end{align}

To evaluate the right hand side, first we observe that the integration away from the two bubbling disks is $O(e^{-\mu_2^k}\delta_k)$ based on the asymptotic behavior of $V_k$. So we focus on the integration over $B(Q_l^k,\epsilon)$ for $l=0,1$.

If $l=0$, $Q_0^k=e_1$ by definition, we show that the integration over $B(e_1,\epsilon)$ is an insignificant error. 

\medskip

For $l=1$, recall that the local maximum points of $V_k$ are at $1$ and $e^{i\pi}$. By standard estimates for single Liouville equations, we have
\begin{align}\label{other-1}
&\frac 1{2\pi}\int_{B(Q_l^k,\epsilon)}(\log \frac{| e_1-\eta|}{|y-\eta|}|\eta|^{2}(\mathfrak{h}_2^k(\delta_k\eta)e^{v_2^k}-\mathfrak{h}_2^k(\delta_ke_1)e^{V_k})d\eta \nonumber\\
=&4\log \frac{|e_1-Q_l^k|}{|y-Q_l^k|}-4\log \frac{|e_1-e^{i\pi}|}{|y-e^{i\pi}|}+O(\mu_2^ke^{-\mu_2^k})\nonumber\\
=&4\log \frac{|1-Q_l^k|}{2}-4\log \frac{|y-Q_l^k|}{|y+1|}+O(\mu_2^ke^{-\mu_2^k}).
\end{align}

Finally for $B(e_1,\tau)$ we claim to have
\begin{equation}\label{april-1}
\frac 1{2\pi}\int_{B(e_1,\tau)}\log |e_1-\eta||\eta|^{2}(\mathfrak{h}_2^k(\delta_k\eta)e^{v_2^k}-\mathfrak{h}_2^k(\delta_ke_1)e^{V_k})d\eta =E
\end{equation}
and for $y\in B_3$ away from bubbling disks,
\begin{equation} \label{april-2}
\frac 1{2\pi}\int_{B(e_1,\tau)}\log |y-\eta||\eta|^{2}(\mathfrak{h}_2^k(\delta_k\eta)e^{v_2^k}-\mathfrak{h}_2^k(\delta_ke_1)e^{V_k})d\eta=E.
\end{equation}

The proof of (\ref{april-2}) follows easily from the same sample computation in the derivation of (\ref{sample-2}). To prove (\ref{april-1}) we first remark that by
scaling around $e_1$ and using the expansion of bubbles for single Liouville equation, we have
$$\int_{B(e_1,\tau)}\log |e_1-\eta ||\eta |^{2}\mathfrak{h}_2^k(\delta_k\eta)e^{v_2^k}d\eta=
-8\pi\mu_2^k+c_k(\mu_2^k)^2e^{-\mu_2^k}+O(\mu_2^ke^{-\mu_2^k})+O(\delta_k^2) $$
where $c_k$ is a constant coming from $\mathfrak{h}_2^k$. After computing the term with $\mathfrak{h}_2^k(\delta_ke_1)$ we see that after cancellation the left-over becomes
$$O(\delta_k(\mu_2^k)^2e^{-\mu_2^k})+O(\mu_2^ke^{-\mu_2^k})+O(\delta_k^2).$$ 
Thus (\ref{april-1}) is verified because by Cauchy's inequality $O(\delta_k (\mu_2^k)^2e^{-\mu_2^k})$ is bounded by the two other terms.
Putting the estimates in different regions together we have
\begin{equation}\label{e-wk-important}
w_k(y)=(4\log \frac{|e_1-Q_1^k|}{2}-4\log \frac{|y-Q_1^k|}{|y+e_1|})+O(\mu_2^ke^{-\mu_2^k})+O(\delta_k^2).
\end{equation}

Now we consider the Pohozaev identity around $Q_1^k$. Let $\Omega_{1,k}=B(Q_1^k,\tau)$ for small $\tau>0$. For $v_2^k$ we have
\begin{align}\label{pi-vk}
&\int_{\Omega_{1,k}}\partial_{\xi}(|y|^{2}\mathfrak{h}_2^k(\delta_ky))e^{v_2^k}-\int_{\partial\Omega_{1,k}}e^{v_2^k}|y|^{2}\mathfrak{h}_2^k(\delta_ky)(\xi\cdot \nu)\\
=&\int_{\partial \Omega_{1,k}}(\partial_{\nu}v_2^k\partial_{\xi}v_2^k-\frac 12 |\nabla v_2^k|^2(\xi \cdot \nu))dS.\nonumber
\end{align}
where $\xi$ is an arbitrary unit vector. Correspondingly the Pohozaev identity for $V_k$ is
\begin{align}\label{pi-Vk}
&\int_{\Omega_{1,k}}\partial_{\xi}(|y|^{2}\mathfrak{h}_2^k(\delta_ke_1))e^{V_k}-\int_{\partial\Omega_{s,k}}e^{V_k}|y|^{2}\mathfrak{h}_2^k(\delta_ke_1)(\xi\cdot \nu)\\
=&\int_{\partial \Omega_{1,k}}(\partial_{\nu}V_k\partial_{\xi}V_k-\frac 12 |\nabla V_k|^2(\xi \cdot \nu))dS.\nonumber
\end{align}

 The second term on the left hand side of (\ref{pi-vk}) or (\ref{pi-Vk}) is $O(e^{-\mu_2^k})$, so both terms are considered as errors. Now we first focus on (\ref{pi-vk}).
  After using the expansion for single equation as in the evaluation of (\ref{pi-e-1}), the first term on the left hand side of (\ref{pi-vk}) is:
\begin{align*}
&\int_{\Omega_{1,k}}\partial_{\xi}(|y|^{2}\mathfrak{h}_2^k(\delta_ky))e^{v_2^k}\\
=&\int_{\Omega_{1,k}}\partial_{\xi}\bigg (\log (|y|^{2}\mathfrak{h}_2^k(\delta_ky)) \bigg )|y|^{2}\mathfrak{h}_2^k(\delta_ky)e^{v_2^k}\\
=&8\pi(2\frac{Q_1^k}{|Q_1^k|^2}+\delta_k\nabla (\log \mathfrak{h}_2^k)(\delta_kQ_s^k))\cdot \xi+E.
\end{align*}
Here we recall that $E=O(\delta_k^2)+O(\mu_2^ke^{-\mu_2^k})$. 
In a similar fashion, the first term of the left hand side of (\ref{pi-Vk}) is
$$\int_{\Omega_{1,k}}\partial_{\xi}(|y|^{2}\mathfrak{h}_2^k(\delta_ke_1))e^{V_k}
=16\pi( e^{i\pi}\cdot \xi)+E.
$$
Recall that  $Q_1^k=-(e_1+m_1^k)+E$, we have
\begin{align*}
&\frac{Q_l^k}{|Q_l^k|^2}+e_1=\frac{1}{\bar Q_l^k}+e_1\\
=&\bar m_1^k+E
=-L\delta_k+E.
\end{align*}
Here we used (\ref{m-2}). Thus the difference of the left hand sides gives the following leading term:
$$-8\pi \delta_kL \cdot \xi+E. $$

To evaluate the right hand side, since $v_2^k=V_k+w_k+E$, the right hand side of (\ref{pi-vk}) is
\begin{align*}
&\int_{\partial\Omega_{s,k}}(\partial_{\nu} v_2^k\partial_{\xi} v_2^k-\frac 12 |\nabla v_2^k|^2(\xi\cdot \nu))dS\\
=&\int_{\partial \Omega_{s,k}}(\partial_{\nu} V_k\partial_{\xi} V_k-\frac 12 |\nabla  V_k|^2(\xi\cdot \nu))dS\\
&+\int_{\partial\Omega_{s,k}}(\partial_{\nu} V_k\partial_{\xi}w_k+\partial_{\nu}w_k\partial_{\xi} V_k-(\nabla  V_k\cdot \nabla w_k)(\xi \cdot \nu))dS+E.
\end{align*}
where we have used $w_k(y)=O(\delta_k)+E$.   Thus the difference of two Pohozaev identities gives
\begin{align}\label{important-1}
&-8\pi \delta_k\partial_{\xi}(\log \mathfrak{h}_2^k)(0)+E\\
=&\int_{\partial\Omega_{1,k}}(\partial_{\nu} V_k\partial_{\xi}w_k+\partial_{\nu}w_k\partial_{\xi} V_k-(\nabla  V_k\cdot \nabla w_k)(\xi \cdot \nu))dS\nonumber
\end{align}

Now we evaluate $\nabla w_k$ on $\partial B(Q_1^k,r)$ for $r>0$ fixed. For simplicity we omit $k$ in $m_1^k$ and $Q_1^k$.  For $y=-e_1+re^{i\theta}$ we obtain
from (\ref{e-wk-important}) that
\begin{align*}
\nabla w_k(y)&=(-4\frac{y-Q_1}{|y-Q_1|^2}+4\frac{y+e_1}{|y+e_1|^2})+E\\
&=4(\frac{1}{\bar y+e_1}-\frac{1}{\bar y-\bar Q_1})+E\\
&=4\frac{-e_1-\bar Q_1}{(\bar y+e_1)(\bar y-\bar Q_1^k)}+E\\
&=4\frac{\bar m_1}{(\bar y+e_1)^2}+E.
\end{align*}
Then for $y=-e_1+re^{i\theta}$, we have
$$
\nabla w_k(y)=\frac{4}{r^2}e^{2i\theta}\bar m_1+E.
$$
If we use $m_1=a+ib$, we can write $\nabla w_k(y)$ as
\begin{equation}\label{grad-w-k}
\nabla w_k(y)=\frac 4{r^2}((a\cos 2\theta+b\sin 2\theta)+i(a\sin 2\theta -b\cos 2\theta))+E.
\end{equation}

On the other hand from the expression of $V_k$ in (\ref{Vk-exp}) we have, for $y=-e_1+re^{i\theta}$,
\begin{equation}\label{grad-V}
\nabla V_k(y)=-4\frac{y+e_1}{|y+e_1|^2}+O(\mu_2^ke^{-\mu_2^k})
=-\frac{4 e^{i\theta}}{r}+O(\mu_2^ke^{-\mu_2^k}).
\end{equation}

Correspondingly we have,
\begin{equation}\label{d-nu-v}
\partial_{\nu}V_k=-\frac{4}{r}+E.
\end{equation}
For $\partial_{\xi}w_k$ we have, using $m_1=a+ib$ in (\ref{grad-w-k})
\begin{equation}\label{d-w-xi}
\partial_{\xi}w_k
=\frac{4}{r^2}\big ((a\cos(2\theta)+b\sin(2\theta))\xi_1
+(a\sin (2\theta)-b\cos(2\theta))\xi_2\big )+E.
\end{equation}

Similarly
\begin{align}\label{d-nu-w}
&\quad \partial_{\nu}w_k \nonumber\\
&=\frac{4}{r^2}[(a\cos(2\theta)+b\sin(2\theta))\cos\theta+(a\sin(2\theta)-b\cos (2\theta) )\sin \theta]
+E \\
&=\frac{4}{r^2}(a\cos \theta +b\sin \theta)+E.\nonumber
\end{align}

Using (\ref{d-nu-v}) and (\ref{d-nu-w}) we have
\begin{equation}\label{one-term-right}
\int_{\partial \Omega_{s,k}}\partial_{\nu}V_k\partial_{\xi}w_k=E. 
\end{equation}
To evaluate other terms on the right hand side, we use (\ref{grad-V}) to have
\begin{equation}\label{d-xi-V}
\partial_{\xi}V_k=-\frac{4}{r}(\cos \theta \xi_1+\sin \theta \xi_2)+E
\end{equation}

By (\ref{d-nu-w}) and (\ref{d-xi-V}) we obtain from direct computation that

\begin{equation}\label{add-3}
\int_{\partial \Omega_{1,k}}\partial_{\nu}w_k\partial_{\xi}V_k=
-\frac{16\pi}{r^2}m_1\cdot \xi.
\end{equation}

Now we compute $\int_{\partial \Omega_{s,k}}(\nabla V_k\cdot \nabla w_k)(\xi\cdot \nu)$.

By the expressions of $\nabla V_k$ and $\nabla w_k$ in (\ref{grad-V}) and (\ref{grad-w-k}), respectively, 
we obtain, after integration that 
\begin{equation}\label{add-4}
\int_{\partial \Omega_{1,k}}(\nabla V_k\cdot \nabla w_k)(\xi\cdot \nu)=-\frac{16\pi}{r^2}m_1\cdot \xi .\nonumber
\end{equation}
Combining (\ref{one-term-right}), (\ref{add-3}) and (\ref{add-4}) we have
\begin{equation}\label{2-part-pi}
\int_{\partial \Omega_{s,k}}(\partial_{\nu}V_k\partial_{\xi}w_k+\partial_{\xi}V_k\partial_{\nu}w_k-(\nabla V_k\cdot \nabla w_k)(\xi\cdot \nu) =E
\end{equation}
Comparing with the left hand side we have
$$\nabla \log \mathfrak{h}_2^k(0)=O(\delta_k^{-1}\mu_2^ke^{-\mu_2^k})+O(\delta_k). $$

\subsection{Improved vanishing rate of the first derivatives}
Our next goal is to prove the following vanishing rate for $\nabla \mathfrak{h}_2^k(0)$:
\begin{equation}\label{vanish-first-tau}
\nabla (\log \mathfrak{h}_2^k)(0)=O(\delta_k\mu_2^k)
\end{equation}

Note that in the previous section we have proved that 
$$\nabla (\log \mathfrak{h}_2^k)(0)=O(\delta_k^{-1}\mu_2^ke^{-\mu_2^k})+O(\delta_k).$$
If $\delta_k\ge C\epsilon_k$, there is nothing to prove. So we assume that 
\begin{equation}\label{delta-small-1}
\delta_k=o(\epsilon_k). 
\end{equation}
By way of contradiction we assume that 
\begin{equation}\label{assum-delta-h}
|\nabla \mathfrak{h}_2^k(0)|/(\delta_k \mu_2^k)\to \infty. 
\end{equation}

Another observation is that based on (\ref{m-2}) and the initial vanishing rate of $\mathfrak{h}_2^k(0)$ we have 
$$\epsilon_k^{-1} | Q_1^k-e^{i\pi}|\le C\epsilon_k^{\epsilon}$$
for some small $\epsilon>0$. Thus $\xi_k$ tends to $U$ after scaling. We need this fact in our argument.

Under the assumption (\ref{delta-small-1}) we claim
\begin{equation}\label{vk-Vk-2}
|v_2^k(-1+\epsilon_ky)-V_k(-1+\epsilon_ky)|\le C\epsilon_k^{\epsilon} (1+|y|)^{-1},\quad 0<|y|<\tau \epsilon_k^{-1}.
\end{equation}
for some small constants $\epsilon>0$ and $\tau>0$ both independent of $k$.  The proof of (\ref{vk-Vk-2}) is very similar to Proposition 3.1 of \cite{wei-zhang-plms} and is omitted.

One major step in the proof of the main theorem is the following estimate:

\begin{prop}\label{key-w8-8} Let $w_k=v_2^k-V_k$, then
$$|w_k(y)|\le C\tilde{\delta_k}, \quad y\in \Omega_k:=B(0,\tau \delta_k^{-1}), $$
where $\tilde{\delta_k}=|\nabla \mathfrak{h}_2^k(0)|\delta_k+\delta_k^2\mu_2^k$.
\end{prop}

\noindent{\bf Proof of Proposition \ref{key-w8-8}:}

Obviously we can assume that $|\nabla \mathfrak{h}_2^k(0)|\delta_k> 2\delta_k^2\mu_2^k$ because otherwise there is nothing to prove. 
Now we recall the equation for $v_2^k$ is (\ref{v2-eq-weak}), which can be written as
$$\Delta v_2^k+\mathfrak{h}_2^k(\delta_ky)e^{v_2^k}=O(\mu_1^ke^{-\mu_1^k}). $$
Since $\mu_1^k$ is very large, the right hand side is at least $O(\delta_k^4)$. 

$v_2^k$ is a constant on $\partial B(0,\tau \delta_k^{-1})$. Moreover $v_2^k(e_1)=\mu_2^k$. Recall that $V_k$ defined in (\ref{Vk-exp}) satisfies
$$\Delta V_k+\mathfrak{h}_2^k(\delta_ke_1)|y|^{2}e^{V_k}=0,\quad \mbox{in}\quad \mathbb R^2, \quad \int_{\mathbb R^2}|y|^{2}e^{V_k}<\infty, $$
$V_k$ has its local maximums at $1$ and $-1$, and $V_k(e_1)=\mu_2^k$.
For $|y|\sim \delta_k^{-1}$,
$$V_k(y)=-\mu_2^k-8\log \delta_k^{-1}+C+O(\delta_k^{2}). $$

Let $\Omega_k=B(0,\tau \delta_k^{-1})$,  we shall derive a precise, point-wise estimate of $w_k$ in $B_3\setminus B(Q_1^k,\tau)$ where $\tau>0$ is a small number independent of $k$. Here we note that among the two local maximum points, we already have $e_1$ as a common local maximum point for both $v_2^k$ and $V_k$ and we shall prove that $w_k$ is very small in $B_3$ if we exclude all bubbling disks except the one around $e_1$. 

Now we write the equation of $w_k$ as
\begin{equation}\label{eq-wk}
\Delta w_k+\mathfrak{h}_2^k(\delta_k y)|y|^{2}e^{\xi_k}w_k=(\mathfrak{h}_2^k(\delta_k e_1)-\mathfrak{h}_2^k(\delta_k y))|y|^{2}e^{V_k}+E_2
\end{equation}
 in $\Omega_k$, where $\xi_k$ is obtained from the mean value theorem:
$$
e^{\xi_k(x)}=\left\{\begin{array}{ll}
\frac{e^{v_2^k(x)}-e^{V_k(x)}}{v_2^k(x)-V_k(x)},\quad \mbox{if}\quad v_2^k(x)\neq V_k(x),\\
\\
e^{V_k(x)},\quad \mbox{if}\quad v_2^k(x)=V_k(x).
\end{array}
\right.
$$
An equivalent form is
\begin{equation}\label{xi-k}
e^{\xi_k(x)}=\int_0^1\frac d{dt}e^{tv_2^k(x)+(1-t)V_k(x)}dt=e^{V_k(x)}\big (1+\frac 12w_k(x)+O(w_k(x)^2)\big ).
\end{equation}
For convenience we write the equation for $w_k$ as
\begin{equation}\label{eq-wk-2}
\Delta w_k+\mathfrak{h}_2^k(\delta_k y)|y|^{2}e^{\xi_k}w_k=\delta_k\nabla \mathfrak{h}_2^k(\delta_k e_1)\cdot (e_1-y)|y|^{2}e^{V_k}+E_2
\end{equation}
where $$E_2=O(\delta_k^2)|y-e_1|^2|y|^{2}e^{V_k}+O(\mu_1^ke^{-\mu_1^k}),\quad y\in \Omega_k. $$
Note that the oscillation of $w_k$ on $\partial \Omega_k$ is $O(\delta_k^{2})$, which all comes from the oscillation of $V_k$. 
It is a trivial and important fact that $E_2$ is small round $e_1$ but greater around $Q_1$. This is the reason that our analysis will first be carried out around $e_1$ and pass to $Q_1$ later. 

Let $M_k=\max_{x\in \bar \Omega_k}|w_k(x)|$. We shall get a contradiction by assuming $M_k/\tilde{\delta_k}\to \infty$. This assumption implies
\begin{equation}\label{big-mk}
M_k/(|\nabla \mathfrak{h}_2^k(0)|\delta_k)\to \infty,\qquad M_k/(\delta_k^2 \mu_2^k)\to \infty. 
\end{equation}
Set
$$\tilde w_k(y)=w_k(y)/M_k,\quad x\in \Omega_k. $$
Clearly $\max_{x\in \Omega_k}|\tilde w_k(x)|=1$. The equation for $\tilde w_k$ is
\begin{equation}\label{t-wk}
\Delta \tilde w_k(y)+|y|^{2}\mathfrak{h}_2^k(\delta_k e_1)e^{\xi_k}\tilde w_k(y)=a_k\cdot (e_1-y)|y|^{2}e^{V_k}+E_3,
\end{equation}
in $\Omega_k$,
where $a_k=\delta_k\nabla \mathfrak{h}_k(0)/M_k\to 0$,
\begin{equation}\label{t-ek}
E_3=o(1)|y-e_1|^2|y|^{2}e^{V_k}+O(\delta_k^2(\mu_2^k)^{-1}),\quad y\in \Omega_k.
\end{equation}
where we used $\mu_1^ke^{-\mu_1^k}/M_k=o(\delta_k^2 (\mu_2^k)^{-1}$. 

Also on the boundary, since $M_k/\tilde \delta_k\to \infty$, we have 
\begin{equation}\label{w-bry}
\tilde w_k=C+o(1/\mu_2^k),\quad \mbox{on}\quad \partial \Omega_k. 
\end{equation}

 By (\ref{vk-Vk-2})
\begin{equation}\label{xi-V-c}
\xi_k(e_1+\epsilon_k z)=V_k(e_1+\epsilon_kz)+O(\epsilon_k^{\epsilon})(1+|z|)^{-1}
\end{equation}

Since $V_k$ is not exactly symmetric around $e_1$, we shall replace the re-scaled version of $V_k$ around $e_1$ by a radial function.
Let $U_k$ be solutions of
\begin{equation}\label{global-to-use}
\Delta U_k+\mathfrak{h}_2^k(\delta_ke_1)e^{U_k}=0,\quad \mbox{in}\quad \mathbb R^2, \quad U_k(0)=\max_{\mathbb R^2}U_k=0.
\end{equation}
By the classification theorem of Caffarelli-Gidas-Spruck \cite{CGS} we have
$$U_k(z)=\log \frac{1}{(1+\frac{\mathfrak{h}_2^k(\delta_ke_1)}{8}|z|^2)^2}$$
and standard refined estimates yield (see \cite{chenlin1,zhangcmp,gluck})
\begin{equation}\label{Vk-rad}
V_k(e_1+\epsilon_k z)+2\log \epsilon_k=U_k(z)+O(\epsilon_k)|z|+O((\mu_2^k)^2\epsilon_k^2).
\end{equation}
Also we observe that
\begin{equation}\label{log-rad}
\log |e_1+\epsilon_k z|=O(\epsilon_k)|z|.
\end{equation}

Thus, the combination of (\ref{xi-V-c}), (\ref{Vk-rad}) and (\ref{log-rad}) gives
\begin{align}\label{xi-U}
&2\log |e_1+\epsilon_kz|+\xi_k(e_1+\epsilon_k z)+2\log \epsilon_k-U_k(z)\\
=&O(\epsilon_k^{\epsilon})(1+|z|)\quad 0\le |z|<\delta_0 \epsilon_k^{-1}.
 \nonumber
\end{align}
for a small $\epsilon>0$ independent of $k$. 
Since we shall use the re-scaled version, based on (\ref{xi-U}) we have
\begin{equation}\label{xi-eU}
\epsilon_k^2 |e_1+\epsilon_k z|^{2}e^{\xi_k(e_1+\epsilon_k z)}
= e^{U_k(z)}+O( \epsilon_k^{\epsilon})(1+|z|)^{-3}
\end{equation}
Here we note that the estimate in (\ref{xi-U}) is not optimal.  In the following we shall put the proof of Proposition \ref{key-w8-8} into a few estimates. In the first estimate we prove

\begin{lem}\label{w-around-e1} For $\delta>0$ small and independent of $k$,
\begin{equation}\label{key-step-1}
\tilde w_k(y)=o(1),\quad \nabla \tilde w_k=o(1) \quad \mbox{in}\quad B(e_1,\delta)\setminus B(e_1,\delta/8)
\end{equation}
where $B(e_1,3\delta)$ does not include other blowup points.
\end{lem}

\noindent{\bf Proof of Lemma  \ref{w-around-e1}:}

If (\ref{key-step-1}) is not true, we have, without loss of generality that $\tilde w_k\to c>0$. This is based on the fact that $\tilde w_k$ tends to a global harmonic function with removable singularity. So $\tilde w_k$ tends to constant. Here we assume $c>0$ but the argument for $c<0$ is the same. Let
\begin{equation}\label{w-ar-e1}
W_k(z)=\tilde w_k(e_1+\epsilon_kz), \quad \epsilon_k=e^{-\frac 12 \mu_2^k},
\end{equation}
then if we use $W$ to denote the limit of $W_k$, we have
$$\Delta W+e^UW=0, \quad \mathbb R^2, \quad |W|\le 1, $$
and $U$ is a solution of $\Delta U+e^U=0$ in $\mathbb R^2$ with $\int_{\mathbb R^2}e^U<\infty$. Since $0$ is the local maximum of $U$,
$$U(z)=\log \frac{1}{(1+\frac 18|z|^2)^2}. $$
Here we further claim that $W\equiv 0$ in $\mathbb R^2$ because $W(0)=|\nabla W(0)|=0$, a fact well known based on the classification of the kernel of the linearized operator. Going back to $W_k$, we have
$$W_k(z)=o(1),\quad |z|\le R_k \mbox{ for some } \quad R_k\to \infty. $$

Based on the expression of $\tilde w_k$, (\ref{Vk-rad}) and (\ref{xi-eU}) we write the equation of $W_k$ as
\begin{equation}\label{e-Wk}
\Delta W_k(z)+\mathfrak{h}_2^k(\delta_ke_1)e^{U_k(z)}W_k(z)=E_2^k,
\end{equation}
for $|z|<\delta_0 \epsilon_k^{-1}$ where a crude estimate of the error term $E_2^k$ is
\begin{equation*}
E_2^k(z)=o(1)\epsilon_k^{\epsilon}(1+|z|)^{-3}.
\end{equation*}

Let
\begin{equation}\label{for-g0}
g_0^k(r)=\frac 1{2\pi}\int_0^{2\pi}W_k(r,\theta)d\theta.
\end{equation}
Then clearly $g_0^k(r)\to c>0$ for $r\sim \epsilon_k^{-1}$.
 The equation for $g_0^k$ is
\begin{align*}
&\frac{d^2}{dr^2}g_0^k(r)+\frac 1r \frac{d}{dr}g_0^k(r)+\mathfrak{h}_2^k(\delta_ke_1)e^{U_k(r)}g_0^k(r)=\tilde E_0^k(r)\\
&g_0^k(0)=\frac{d}{dr}g_0^k(0)=0.
\end{align*}
where $\tilde E_0^k(r)$ has the same upper bound as that of $E_2^k(r)$:
$$|\tilde E_0^k(r)|\le o(1)\epsilon_k^{\epsilon}(1+r)^{-3}. $$

For the homogeneous equation, the two fundamental solutions are known: $g_{01}$, $g_{02}$, where
$$g_{01}=\frac{1-c_1r^2}{1+c_1r^2},\quad c_1=\frac{\mathfrak{h}_k(\delta_ke_1)}8.$$
By the standard reduction of order process, $g_{02}(r)=O(\log r)$ for $r>1$.
Then it is easy to obtain, assuming $|W_k(z)|\le 1$, that
\begin{align*}
|g_0(r)|\le C|g_{01}(r)|\int_0^r s|\tilde E_0^k(s) g_{02}(s)|ds+C|g_{02}(r)|\int_0^r s|g_{01}(s)\tilde E_0^k(s)|ds\\
\le C\epsilon_k^{\epsilon}\log (2+r). \quad 0<r<\delta_0 \epsilon_k^{-1}.
\end{align*}
Clearly this is a contradiction to (\ref{for-g0}). We have proved $c=0$, which means $\tilde w_k=o(1)$ in $B(e_1, \delta_0)\setminus B(e_1, \delta_0/8)$.
Then it is easy to use the equation for $\tilde w_k$ and standard Harnack inequality to prove
$\nabla \tilde w_k=o(1)$ in the same region.
Lemma \ref{w-around-e1} is established. $\Box$

\medskip

The second estimate is a more precise description of $\tilde w_k$ around $e_1$:
\begin{lem}\label{t-w-1-better} For any given $\sigma\in (0,1)$ there exists $C>0$ such that
\begin{equation}\label{for-lambda-k}
|\tilde w_k(e_1+\epsilon_kz)|\le C\epsilon_k^{\sigma} (1+|z|)^{\sigma},\quad 0<|z|<\tau \epsilon_k^{-1}.
\end{equation}
for some $\tau>0$.
\end{lem}
\begin{rem}
The reason that we have $\sigma\in (0,1)$ is because we currently don't have a very good estimate of $\tilde w_k$ outside the bubbling area, except knowing it is $o(1)$. Once a better estimate of $\tilde w_k$ outside the bubbling disks is known, (\ref{for-lambda-k})  can be improved. 
\end{rem}

\noindent{\bf Proof of Lemma \ref{t-w-1-better}:} Let $W_k$ be defined as in (\ref{w-ar-e1}). In order to obtain a better estimate we need to write the equation of $W_k$ more precisely than (\ref{e-Wk}):
\begin{equation}\label{w-more}
\Delta W_k+\mathfrak{h}_2^k(\delta_ke_1)e^{\Theta_k}W_k=E_3^k(z), \quad z\in \Omega_{Wk}
\end{equation}
where
$\Theta_k$ is defined by
$$e^{\Theta_k(z)}=|e_1+\epsilon_k z|^{2}e^{\xi_k(e_1+\epsilon_kz)+2\log \epsilon_k}, $$
$\Omega_{Wk}=B(0,\tau \epsilon_k^{-1})$ and $E_3^k(z)$ satisfies
$$E_3^k(z)=O(\epsilon_k)(1+|z|)^{-3},\quad z\in \Omega_{Wk}. $$
Here we observe that by Lemma \ref{w-around-e1}  $W_k=o(1)$ on $\partial \Omega_{Wk}$. 
Let
$$\Lambda_k=\max_{z\in \Omega_{Wk}}\frac{|W_k(z)|}{\epsilon_k^{\sigma}(1+|z|)^{\sigma}}. $$
If (\ref{for-lambda-k}) does not hold, $\Lambda_k\to \infty$ and we use $z_k$ to denote where $\Lambda_k$ is attained. Note that because of the smallness of $W_k$ on $\partial \Omega_{Wk}$, $z_k$ is an interior point. Let
$$g_k(z)=\frac{W_k(z)}{\Lambda_k (1+|z_k|)^{\sigma}\epsilon_k^{\sigma}},\quad z\in \Omega_{Wk}, $$
we see immediately that
\begin{equation}\label{g-sub-linear}
|g_k(z)|=\frac{|W_k(z)|}{\epsilon_k^{\sigma}\Lambda_k(1+|z|)^{\sigma}}\cdot \frac{(1+|z|)^{\sigma}}{(1+|z_k|)^{\sigma}}\le  \frac{(1+|z|)^{\sigma}}{(1+|z_k|)^{\sigma}}.
\end{equation}
Note that $\sigma$ can be as close to $1$ as needed. The equation of $g_k$ is
$$\Delta g_k(z)+\mathfrak{h}_2^k(\delta_k e_1)e^{\Theta_k}g_k=o(\epsilon_k^{1-\sigma})\frac{(1+|z|)^{-3}}{(1+|z_k|)^{\sigma}}, \quad \mbox{in}\quad \Omega_{Wk}. $$
Then we can obtain a contradiction to $|g_k(z_k)|=1$ as follows: If $\lim_{k\to \infty}z_k=P\in \mathbb R^2$, this is not possible because that fact that $g_k(0)=|\nabla g_k(0)|=0$ and the sub-linear growth of $g_k$ in (\ref{g-sub-linear}) implies that $g_k\to 0$ over any compact subset of $\mathbb R^2$ (see \cite{chenlin1,zhangcmp}). So we have $|z_k|\to \infty$. But this would lead to a contradiction again by using the Green's representation of $g_k$:
\begin{align} \label{temp-1}
&\pm 1=g_k(z_k)=g_k(z_k)-g_k(0)\\
&=\int_{\Omega_{k,1}}(G_k(z_k,\eta)-G_k(0,\eta))(\mathfrak{h}_2^k(\delta_k e_1)e^{\Theta_k}g_k(\eta)+o(\epsilon_k^{1-\sigma})\frac{(1+|\eta |)^{-3}}{(1+|z_k|)^{\sigma}})d\eta+o(1).\nonumber
\end{align}
where $G_k(y,\eta)$ is the Green's function on $\Omega_{Wk}$ and $o(1)$ in the equation above comes from the smallness of $W_k$ on $\partial \Omega_{Wk}$. Let $L_k=\tau\epsilon_k^{-1}$, the expression of $G_k$ is 
$$G_k(y,\eta)=-\frac{1}{2\pi}\log |y-\eta|+\frac 1{2\pi}\log (\frac{|\eta |}{L_k}|\frac{L_k^2\eta}{|\eta |^2}-y|). $$
$$G_k(z_k,\eta)-G_k(0,\eta)=-\frac{1}{2\pi}\log |z_k-\eta |+\frac 1{2\pi}\log |\frac{z_k}{|z_k|}-\frac{\eta z_k}{L_k^2}|+\frac 1{2\pi}\log |\eta |. $$
Using this expression in (\ref{temp-1}) we obtain from elementary computation that the right hand side of (\ref{temp-1}) is $o(1)$, a contradiction to $|g_k(z_k)|=1$. Lemma \ref{t-w-1-better} is
established. $\Box$

\medskip

The smallness of $\tilde w_k$ around $e_1$ can be used to obtain the following third key estimate:
\begin{lem}\label{small-other}
\begin{equation}\label{key-step-2}
\tilde w_k=o(1)\quad \mbox{in}\quad B(e^{i\pi},\tau).
\end{equation}
\end{lem}

\noindent{\bf Proof of Lemma \ref{small-other}:}
We abuse the notation $W_k$ by defining it as
$$W_k(z)=\tilde w_k(e^{i\pi}+\epsilon_k z),\quad z\in B(0,\tau \epsilon_k^{-1}). $$
Here we point out that based on (\ref{m-2}) and (\ref{delta-small-1}) we have $\epsilon_k^{-1}|Q_1^k-e^{i\pi}|\to 0$. So the scaling around $e^{i\pi}$ or $Q_1^k$ does not affect the
limit function.

$$\epsilon_k^2 |e^{i\pi}+\epsilon_kz|^{2}\mathfrak{h}_2^k(\delta_ke_1)e^{\xi_k(e^{i\pi}+\epsilon_kz)}\to e^{U(z)} $$
where $U(z)$ is a solution of
$$\Delta U+e^U=0,\quad \mbox{in}\quad \mathbb R^2, \quad \int_{\mathbb R^2}e^U<\infty. $$
Here we recall that $\lim_{k\to \infty} \mathfrak{h}_2^k(\delta_k e_1)=1$.
Since $W_k$ converges to a solution of the linearized equation:
$$\Delta W+e^UW=0, \quad \mbox{in}\quad \mathbb R^2. $$
$W$ can be written as a linear combination of three functions:
$$W(x)=c_0\phi_0+c_1\phi_1+c_2\phi_2, $$
where
$$\phi_0=\frac{1-\frac 18 |x|^2}{1+\frac 18 |x|^2} $$
$$\phi_1=\frac{x_1}{1+\frac 18 |x|^2},\quad \phi_2=\frac{x_2}{1+\frac 18|x|^2}. $$
The remaining part of the proof consists of proving $c_0=0$ and $c_1=c_2=0$. First we prove $c_0=0$.

\noindent{\bf Step one: $c_0=0$.}
First we write the equation for $W_k$ in a convenient form. Since
$$|e^{i\pi}+\epsilon_kz|^{2}\mathfrak{h}_2^k(\delta_ke_1)=\mathfrak{h}_2^k(\delta_ke_1)+O(\epsilon_k z),$$
and
$$\epsilon_k^2e^{\xi_k(e^{i\pi}+\epsilon_kz)}=e^{U_k(z)}+O(\epsilon_k^{\epsilon})(1+|z|)^{-3}. $$
Based on (\ref{t-wk}) we write the equation for $W_k$ as
\begin{equation}\label{around-l}
\Delta W_k(z)+\mathfrak{h}_2^k(\delta_ke_1)e^{U_k}W_k=E_l^k(z)
\end{equation}
where
$$E_l^k(z)=O(\epsilon_k^{\epsilon})(1+|z|)^{-3}\quad \mbox{in}\quad \Omega_{k,l}.$$
In order to prove $c_0=0$, the key is to control the derivative of $W_0^k(r)$ where
$$W_0^k(r)=\frac 1{2\pi r}\int_{\partial B_r} W_k(re^{i\theta})dS, \quad 0<r<\tau \epsilon_k^{-1}. $$
To obtain a control of $\frac{d}{dr}W_0^k(r)$ we use $\phi_0^k(r)$ as the radial solution of
$$\Delta \phi_0^k+\mathfrak{h}_2^k(\delta_k e_1)e^{U_k}\phi_0^k=0, \quad \mbox{in }\quad \mathbb R^2. $$

When $k\to \infty$, $\phi_0^k\to c_0\phi_0$. Thus using the equation for $\phi_0^k$ and $W_k$, we have
\begin{equation}\label{c-0-pf}
\int_{\partial B_r}(\partial_{\nu}W_k\phi_0^k-\partial_{\nu}\phi_0^kW_k)=o(\epsilon_k^{\epsilon}). \end{equation}

Thus from (\ref{c-0-pf}) we have
\begin{equation}\label{W-0-d}
\frac{d}{dr}W_0^k(r)=\frac{1}{2\pi r}\int_{\partial B_r}\partial_{\nu}W_k=o(\epsilon_k^{\epsilon})/r+O(1/r^3),\quad 1<r<\tau \epsilon_k^{-1}.
\end{equation}
Since we have known that
$$W_0^k(\tau \epsilon_k^{-1})=o(1). $$
By the fundamental theorem of calculus we have
$$W_0^k(r)=W_0^k(\tau\epsilon_k^{-1})+\int_{\tau \epsilon_k^{-1}}^r(\frac{o(\epsilon_k^{\epsilon})}{s}+O(s^{-3}))ds=O(1/r^2)+O(\epsilon_k^{\epsilon}\log
\frac{1}{\epsilon_k}) $$
for $r\ge 1$. Thus
$c_0=0$ because $W_0^k(r)\to c_0\phi_0$, which means when $r$ is large, it is $-c_0+O(1/r^2)$.

\medskip

\noindent{\bf Step two $c_1=c_2=0$}.
 We first observe that Lemma \ref{small-other} follows from this. Indeed, once we have proved $c_1=c_2=c_0=0$ around $1$ and $-1$, it is easy to use maximum principle to prove $\tilde w_k=o(1)$ in $B_3$ using $\tilde w_k=o(1)$ on $\partial B_3$ and the Green's representation of $\tilde w_k$. The smallness of $\tilde w_k$ immediately implies $\tilde w_k=o(1)$ in $B_R$ for any fixed $R>>1$. Outside $B_R$, a crude estimate of $v_k$ is
  $$v_2^k(y)\le -\mu_2^k-8\log |y|+C, \quad 3<|y|<\tau \delta_k^{-1}. $$
  Using this and the Green's representation of $w_k$ we can first observe that the oscillation on each $\partial B_r$ is $o(1)$ ($R<r<\tau \delta_k^{-1}/2$) and then by the Green's representation of $\tilde w_k$ and fast decay rate of $e^{V_k}$ we obtain $\tilde w_k=o(1)$ in $\overline{B(0,\tau \delta_k^{-1})}$. A contradiction to $\max |\tilde w_k|=1$.

 There are 2 local maximums with one of them being $e_1$. Correspondingly there are 2 global solutions $V_k$ and $\bar V_k$ that
approximate $v_2^k$ accurately near $Q_0^k=e_1$ and $Q_1^k$, respectively. For $\bar V_k$ the expression is
\begin{equation}\label{bar-Vk}
\bar V_k=\log \frac{e^{\bar \mu_2^k}}{(1+\frac{e^{\bar \mu_2^k}}{\bar D_k}|y^{2}-(e_1+p_k)|^2)^2},\quad l=0,...,N, 
\end{equation}
where $p_k=E$ and
\begin{equation}\label{def-D}
D_k=\frac{32}{\mathfrak{h}_2^k(\delta_ke_1)},\quad \bar D_k=\frac{32}{\mathfrak{h}_2^k(\delta_kQ_1^k)}.
\end{equation}
The equation that $\bar V_{k}$ satisfies is 
$$\Delta \bar V_{k}+|y|^{2}\mathfrak{h}_2^k(\delta_k Q_1^k)e^{\bar V_{k}}=0,\quad \mbox{in}\quad \mathbb R^2. $$
Since $v_2^k$ and $\bar V_{k}$ have the same common local maximum at $Q_1^k$, it is easy to see from the expression of $\bar V_k$ in (\ref{bar-Vk}) that
\begin{equation}\label{ql-exp}
Q_1^k=-1-\frac{p_k}{2}+O(|p_k|^2),
\end{equation}
Let $\bar M_{k}$ be the maximum of $|v_2^k-\bar V_{k}|$ and we claim that $\bar M_{k}$ and $M_k$ are comparable:
\begin{equation}\label{M-comp}
\bar M_{k}\sim M_{k},
\end{equation}

The proof of (\ref{M-comp}) is as follows: We use $\bar L$ to denote the limit of $(v_2^k-\bar V_{k})/\bar M_{k}$ around $e_1$:
$$\frac{(v_2^k-\bar V_{k})(e_1+\epsilon_kz)}{\bar M_{k}}=\bar L+o(1),\quad |z|\le \tau \epsilon_k^{-1} $$
where
$$ \bar L=\bar c_{1}\frac{z_1}{1+\frac 18 |z|^2}+\bar c_{2}\frac{z_2}{1+\frac 18 |z|^2}$$
If both $\bar c_{1}$ and $\bar c_{2}$ are zero, we can obtain a contradiction just like the beginning of step two. So at least one of them is not zero.
By Lemma \ref{t-w-1-better} we have
\begin{equation}\label{Q-bad}
v_2^k(e_1+\epsilon_kz)-V_{k}(e_1+\epsilon_kz)=O(\epsilon_k^{\sigma})(1+|z|)^{\sigma} M_{k},\quad |z|<\tau \epsilon_k^{-1}.
\end{equation}

To determine $\bar L$ we see that
\begin{align*}
    &\frac{v_2^k(e_1+\epsilon_k z)-\bar V_{k}(e_1+\epsilon_kz)}{\bar M_k}\\
    =&o(\epsilon_k^{\sigma})(1+|z|)^{\sigma}+\frac{V_{k}(e_1+\epsilon_kz)-\bar V_{k}(e_1+\epsilon_kz)}{\bar M_k}. 
\end{align*}
This expression says that $\bar L$ is mainly determined by the difference of two global solutions $V_{k}$ and $\bar V_{k}$. In order to obtain a contradiction to our assumption we will put the difference in several terms. The main idea in this part of the reasoning is that ``first order terms" tell us what the kernel functions should be, then the ``second order terms" tell us where the pathology is. 

We write $V_{k}(y)-\bar V_{k}(y)$ as
$$V_{k}(y)-\bar V_{k}(y)=\mu_2^k-\bar \mu_2^k+2A-A^2+O(|A|^3) $$
where
$$A(y)=\frac{\frac{e^{\bar \mu_2^k}}{\bar D_k}|y^{2}-e_1-p_k|^2-\frac{e^{\mu_2^k}}{D_k}|y^{2}-e_1|^2}{1+\frac{e^{\mu_2^k}}{D_k}|y^{2}-e_1|^2}.$$
Here for convenience we abuse the notation $\epsilon_k$ by assuming $\epsilon_k=e^{-\mu_2^k/2}$. From $A$ we claim that

\begin{align}\label{late-1}
&V_{k}(e_1+\epsilon_kz)-\bar V_{k}(e_1+\epsilon_kz)\\
=&\phi_1+\phi_2+\phi_3+\phi_4+\mathfrak{R},\nonumber
\end{align}
where
\begin{align*}
&\phi_1=(\mu_2^k-\bar \mu_2^k)(1-\frac{4}{D_k}|z+O(\epsilon_k)|z|^2|^2)/B, \\
&\phi_2=\frac{8}{D_k} \delta_k\nabla\mathfrak{h}_2^k(\delta_k e_1)(Q_1^k-e_1)|z|^2/B\\
&\phi_3=\frac{8}{D_k B}Re((z+O(\epsilon_k|z|^2)(\frac{-\bar p_k}{\epsilon_k}))
\\
&\phi_4=\frac{|p_k|^2}{\epsilon_k^2}\bigg (\frac{2}{D_k  B}-\frac{8|z|^2}{D_k^2 B^2}
-\frac{8}{D_k^2B^2}|z|^2\cos (2\theta-2\theta_{1}^k) \bigg ),\\
&B=1+\frac{4}{D_k}|z+O(\epsilon_k|z|^2)|^2,
\end{align*}
and $\mathfrak{R}_k$ is the collections of other insignificant terms. $\theta_1^k$ and $\theta$ come from
$$|p_k|e^{\theta_1^k}=p_k,\quad z=|z|e^{i\theta}.$$
The proof of (\ref{late-1}) can be found in Appendix A. 
Here $\phi_1$, $\phi_3$ correspond to solutions to the linearized operator. Here we note that if we set $\bar \epsilon_{k}=e^{-\bar\mu_2^k/2}$, there is no essential difference in the proof. If $|\mu_2^k-\bar \mu_2^k|/\bar M_k\ge C$ there is no way to obtain a limit in the form of $\bar L$ mentioned before. Thus we must have $|\mu_2^k-\bar \mu_2^k|/\bar M_k\to 0$. After simplification (see $\phi_3$ of (\ref{late-1})) we have
\begin{align}\label{c-12}
\bar c_{1}=\lim_{k\to \infty}\frac{|p_k|}{4\bar M_k\epsilon_k}\cos(\theta_{1}^k),\\
\bar c_{2}=\lim_{k\to \infty}
\frac{|p_k|}{4\bar M_k\epsilon_k}\sin(\theta_{1}^k)\nonumber
\end{align}

Since $|p_k|=E$, an upper bound for $\bar M_k$ is 
\begin{equation}\label{small-M}
\bar M_k\le C\mu_2^k\epsilon_k+C\delta_k^2\epsilon_k^{-1}.
\end{equation}

Equation (\ref{c-12}) gives us a key observation: $|\bar c_{1}|+|\bar c_{2}|\sim |p_k|/(\epsilon_k \bar M_k)$. Recall that
$$\tilde w_k=(v_2^k-V_k)/M_k, $$
$\hat w_k$ would converge to $\displaystyle{\frac{c_1z_1+c_2z_2}{1+\frac 18|z|^2}}$ near $Q_1^k$, with $M_k\sim \frac{|p_k|}{\epsilon_k}$. Then we see that $\bar M_k\sim M_k$. 
 (\ref{M-comp}) is established.  From now on for convenience we shall just use $M_k$.

\medskip

Set 
$$\bar w_{k}=(v_2^k-\bar V_{k})/M_k,$$
then we have $\bar w_{k}(Q_1^k)=|\nabla \bar w_{k}(Q_1^k)|=0$.

The equation of $\bar w_{k}$ can be written as 
\begin{align}\label{wlk-bs}
   &\Delta \bar w_{k}+|y|^{2}\mathfrak{h}_2^k(\delta_k Q_1^k)e^{\xi_1}\bar w_{k}\\
=&\sigma_k\nabla\mathfrak{h}_2^k(\delta_kQ_1^k)(Q_1^k-y)|y|^{2}e^{\bar V_{k}}-\frac{\delta_k^2}{M_k}\sum_{|\alpha |=2}\frac{\partial^{\alpha}\mathfrak{h}_2^k(\delta_kQ_1^k)}{\alpha !}(y-Q_1^k)^{\alpha}
|y|^{2}e^{\bar V_k}\nonumber \\
&+O(\frac{\delta_k^3}{M_k}|y|^2|y-Q_1^k|^2e^{\bar V_k}.  \nonumber
\end{align}
where $\sigma_k=\delta_k/M_k$ and we omitted $k$ in $\xi_1$. $\xi_1$ comes from the Mean Value Theorem and satisfies
\begin{equation}\label{around-ls}
e^{\xi_1}=e^{\bar V_{k}}(1+\frac 12M_k\bar w_{k}+O(M_k^2\bar w_{k}^2)). 
\end{equation}
The function $\bar w_{k}$ satisfies
\begin{equation}\label{aroud-s-1}
\lim_{k\to \infty}\bar w_{k}(e_1+\epsilon_k z)=\frac{\bar c_{1}z_1+\bar c_{2}z_2}{1+\frac 18 |z|^2}.
\end{equation}

First we have
\begin{align}\label{w-precise}
    \bar w_{k}(y)&=\int_{\Omega_k}(G_k(y,\eta)-G_k(Q_1,\eta))(\mathfrak{h}_2^k(\delta_k Q_1)|\eta |^{2}e^{\xi_1}\bar w_{k}(\eta)\\
    &+\sigma_k \nabla \mathfrak{h}_2^k(\delta_k Q_1)(\eta-Q_1)|\eta |^{2}e^{\bar V_{k}}\nonumber\\
&-\frac{\delta_k^2}{M_k}\sum_{|\alpha |=2}\frac{\partial^{\alpha}\mathfrak{h}_2^k(\delta_kQ_1^k)}{\alpha !}(y-Q_1^k)^{\alpha})d\eta+O(\delta_k^{3}/M_k)\nonumber
\end{align}
with $\sigma_k$ satisfying $\sigma_k|\nabla \mathfrak{h}_2^k(0)|=o(1)$.
The evaluation of (\ref{w-precise}) is involved with integration on the two major disks: One around $Q_1$ and the other one around $e_1$. 
For the integration around $e_1$, it is more convenient to replace
$e^{\bar V_k}$ by $e^{V_k}$ based on the following estimate:
$$\frac{\bar V_k(e_1+\epsilon_kz)-V_k(e_1+\epsilon_kz)}{M_k}=\bar w_k(e_1+\epsilon_kz)+o(\epsilon_k^{\sigma}(1+|z|)^{\sigma}). $$
The following identity, which is due to (\ref{late-1}), provides a key computation of (\ref{w-precise}):
\begin{align}\label{late-2-r}
&\int_{B(e_1,\tau)}\mathfrak{h}_2^k(\delta_k Q_1)|\eta |^{2}e^{\xi_1}\bar w_{k}(\eta)
+\sigma_k \nabla \mathfrak{h}_2^k(\delta_k Q_1)(\eta-Q_1)|\eta |^{2}e^{\bar V_{k}}\\
&\qquad -\frac{\delta_k^2}{M_k}\sum_{|\alpha |=2}\frac{\partial^{\alpha}\mathfrak{h}_2^k(\delta_k Q_1^k)}{\alpha !}(\eta-Q_1^k)^{\alpha}|\eta |^2e^{\bar V_k}d\eta \nonumber\\
=&\frac{\pi}{4}\frac{|p_k|^2}{\epsilon_k^2 M_{k}^2}M_k-2\pi\frac{\delta_k^2}{M_k}\frac{\Delta \mathfrak{h}_2^k(\delta_k Q_1^k)}{\mathfrak{h}_2^k(\delta_k Q_1^k)}+o(\epsilon_k).\nonumber
\end{align}
The proof of (\ref{late-2-r}) can be found in Appendix B.

Equation (\ref{late-2-r}) also leads to a more accurate estimate of $\bar w_k$ in regions between bubbling disks. By the Green's representation formula of $\bar w_k$ it is easy to have
\begin{align*}
& \bar w_k(y)
=-\frac{1}{2\pi}\int_{\Omega_k}\log \frac{|y-\eta |}{|Q_1^k-\eta |}\bigg ( \bar w_{k}(\eta )\mathfrak{h}_2^k(\delta_k Q_1^k)|\eta |^{2}e^{\xi_1}\\
& +\sigma_k\nabla \mathfrak{h}_k(\delta_k Q_l^k)(\eta-Q_l^k)
|\eta |^{2}e^{\bar V_{k}}+O(\frac{\delta_k^2}{M_k})|\eta -Q_1^k|^2|\eta |^2e^{\bar V_k}\bigg )d\eta+o(\epsilon_k)
\end{align*}
for $|y|\sim 1$ and $y\not \in (B(Q_1,\tau)\cup B(e_1,\tau))$. 
Writing the logarithmic term as
$$\log \frac{|y-\eta |}{|Q_1^k-\eta |}=\log \frac{|y-1|}{|Q_1^k-1|}+(\log \frac{|y-\eta |}{|Q_1^k-\eta |}-\log \frac{|y-1|}{|Q_1^k-1|}),$$
we see that the integration related to the second term is $O(\epsilon_k)$. The integration involving the first term is $O(\epsilon_k^{\sigma})$ by (\ref{late-2-r}).
Therefore
$$|\bar w_{k}(y)|=o(\frac{1}{\mu_2^k}),\quad y\in B_3\setminus \cup_{s=0}^1B(Q_s^k,\tau). $$
Thus this extra control of $\bar w_{k}$ away from bubbling disks gives a better estimate than (\ref{Q-bad}) around $Q_1^k$: Using the same argument for Lemma \ref{t-w-1-better} we have
\begin{equation}\label{much-better-s}
|\bar w_{k}(Q_1^k+\epsilon_k z)|\le o(\epsilon_k)\frac{(1+|z|)}{\log (2+|z|)},\quad |z|<\tau \epsilon_k^{-1}.
\end{equation}

Now we deduce a better estimate of $\bar w_k$ away from the two bubbling disks. 
\begin{align*}
   & \bar w_k(y)=\int_{\Omega_k} (G_k(y,\eta)-G_k(Q_1^k,\eta))\bigg (\mathfrak{h}_2^k(\delta_k Q_1^k)|\eta |^2e^{\xi_1}\bar w_k(\eta)\\
   & +\sigma_k
    \nabla\mathfrak{h}_2^k(\delta_k Q_1^k)(\eta-Q_1^k)|\eta |^2e^{\bar v_k}
    -\frac{\delta_k^2}{M_k}\sum_{|\alpha |=2}\frac{\partial^{\alpha}\mathfrak{h}_2^k(\delta_kQ_1^k)}{\alpha !}(\eta-Q_1^k)^{\alpha}|\eta|^2e^{\bar V_k}\bigg )d\eta\\
    &+o(\epsilon_k). 
\end{align*}
Note that the boundary term is $O(\delta_k^3)/M_k$. Since $M_k/\delta_k^2\to \infty$ and $\delta_k=o(\epsilon_k)$ we have $O(\delta_k^3/M_k)=o(\epsilon_k)$. 

The integration of the first part is $o(\epsilon_k)$. The integration on the region away from the two bubbling disks is also clearly $o(\epsilon_k)$. So the main focus is the integration on $B(e_1,\tau)$. For this integration we set
$$H_{2,y}(\eta)=G_k(y,\eta)-G_k(Q_1^k,\eta)$$
then the integration over $B(e_1,\tau)$ is now 
\begin{align*}
&H_{2,y}(e_1)\int_{B(e_1,\tau)}
|\eta |^2(\mathfrak{h}_2^k(\delta_kQ_1^k)e^{\xi_1}\bar w_k+\sigma_k\nabla\mathfrak{h}_2^k(\delta_kQ_1^k)(\eta-Q_1^k)e^{\bar V_k}\\
&-\frac{\delta_k^2}{M_k}\sum_{|\alpha |=2}\frac{\partial^{\alpha}\mathfrak{h}_2^k(\delta_kQ_1^k)}{\alpha !}(\eta-Q_1^k)^{\alpha}|\eta|^2e^{\bar V_k})d\eta\\
&+\int_{B(e_1,\tau)}(H_{2,y}(\eta)-H_{2,y}(e_1)|\eta |^2(\mathfrak{h}_2^k(\delta_kQ_1^k)e^{\xi_1}\bar w_k\\
&+\int_{B(e_1,\tau)}(H_{2,y}(\eta)-H_{2,y}(e_1))
 |\eta |^2e^{\bar V_k}\bigg (\sigma_k\nabla\mathfrak{h}_2^k(\delta_kQ_1^k)(\eta-Q_1^k) +\frac{O(\delta_k^2)}{M_k}|\eta-Q_1^k|^2 \bigg ).
\end{align*}
The first term is a fixed quantity determined by (\ref{late-2-r}). The second term is only involved with the integration of 
the term with $\bar w_k$. The other two terms lead to $o(\epsilon_k)$. 
The result of this computation is
$$H_{2,y}(e_1)(\frac{\pi}4(\frac{|p_k|}{\epsilon_k M_k})^2M_k-2\pi\frac{\delta_k^2}{M_k}\frac{\Delta \mathfrak{h}_2^k(0)}{\mathfrak{h}_2^k(0)})+16\pi \epsilon_k(\partial_1H_{2,y}(e_1) \bar c_1+\partial_2 H_{2,y}(e_1)\bar c_2).  $$
Then we estimate 
$$v_2^k(y)=\bar V_k(y)+M_k\bar w_k(y), $$
for $y$ away from two balls:
\begin{align}\label{compare-1}
v_2^k(y)=\bar V_k(y)+M_k\bigg (H_{2,y}(e_1)\frac{\pi}4(\frac{p_k}{\epsilon_k M_k})^2M_k
\\
+16\pi \epsilon_k(\partial_1H_{2,y}(e_1)\bar c_1+\partial_2 H_{2,y}(e_1)\bar c_2)\bigg )
+o(M_k\epsilon_k). \nonumber
\end{align}
Obviously another evaluation is
$$v_2^k=V_k+M_k\tilde w_k. $$
The comparison of these two expressions will lead to a contradiction. 
Let 
$$\hat w_k=\frac{v_2^k-V_k}{M_k}, $$
then by exactly the same computation, let 
$$H_{1,y}(\eta)=G_k(y,\eta)-G_k(e_1,\eta),$$
we have
\begin{align*}
\hat w_k(y)&=H_{1,y}(Q_1^k)\bigg (\frac{\pi}4(\frac{p_k}{\epsilon_kM_k})^2M_k
-2\pi\frac{\delta_k^2}{M_k}\frac{\Delta \mathfrak{h}_2^k(0)}{\mathfrak{h}_2^k(0)}\bigg )\\
&+16\pi \epsilon_k(\partial_1 H_{1,y}(Q_1^k)\bar c_1+\partial_2 H_{1,y}(Q_1^k)\bar c_2)+o(\epsilon_k). 
\end{align*}

Now we determine the limit function of $\frac{v_2^k-V_k}{M_k}$ near $Q_1^k$ based on $(\bar V_k-V_k)/M_k$. Since the difference between $\bar \mu_2^k$ and $\mu_2^k$ is insignificant we use simplified notations: 
$$\bar V_k=\mu_2^k-2\log (1+\frac{e^{\mu_2^k}}{D}|y^2-e_1-p|^2) $$
$$V_k=\mu_2^k-2\log (1+\frac{e^{\mu_2^k}}{D}|y^2-e_1|^2). $$
\begin{align*}
    \bar V_k-V_k=2\log (1+\frac{\frac{e^{\mu_2^k}}{D}(|y^2-e_1|^2-|y^2-e_1-p|^2)}{1+\frac{e^{\mu_2^k}}{D}|y^2-e_1-p|^2})\\
    =2\log (1+A)=2A-A^2+\mbox{insignificant error}.
\end{align*}
where 
$$A=\frac{\frac{e^{\mu_2^k}}{D}(|y^2-e_1|^2-|y^2-e_1-p|^2)}{1+\frac{e^{\mu_2^k}}{D}|y^2-e_1-p|^2}. $$
$y=-1-\frac{p_k}2+\epsilon_kz$, 
$$y^2=1+p_k-2\epsilon_k z+\epsilon_k^2z^2+O(\mu_k\epsilon_k^3z). $$
Then
\begin{align*}
    |y^2-1|^2=4\epsilon_k^2\bigg (|z-\frac{\epsilon_k z^2}{2}|^2-2Re((z-\frac{\epsilon_k z^2}2)\frac{\bar p_k}{2\epsilon_k})+|\frac{p_k}{2\epsilon_k^2}|^2\\
    +O(\mu_k \epsilon_k^2|z|^2)\bigg )
\end{align*}
$$|y^2-1-p_k|^2=4\epsilon_k^2(|z-\frac{\epsilon_k}2z^2|^2+O(\mu_k\epsilon_k^2|z|^2). $$
Using these expressions we have

$$A=\frac{-\frac 18 \bigg (2 Re(z\frac{\bar p_k}{2\epsilon_k})+|\frac{p_k}{2\epsilon_k}|^2+O(\mu_k\epsilon_k^2|z|^2)\bigg )}{1+\frac 18|z|^2}$$

The leading term of $A^2$ is 
$$\frac{|z|^2|\frac{p_k}{2\epsilon_k}|^2\cos^2(\theta-\theta_1)}{64(1+\frac 18 |z|^2)^2}$$

Thus the leading term of $2A-A^2$ is 
$$
\frac{-\frac 14Re(z\frac{\bar p_k}{\epsilon_k})+2|\frac{p_k}{2\epsilon_k}|^2}{1+\frac 18|z|^2}
-\frac{|z|^2|\frac{p_k}{2\epsilon_k}|^2\cos^2(\theta-\theta_1)}{64(1+\frac 18 |z|^2)^2}
+O(\mu_k\epsilon_k^2|z|^2(1+|z|)^{-4}. 
$$
It is also important to observe that 
$$c_1=-\frac 14\frac{|p_k|}{M_k\epsilon_k}\cos \theta_1, \quad c_2=-\frac 14\frac{|p_k|}{M_k\epsilon_k}\sin \theta_1. $$

Corresponding to (\ref{compare-1}) we also have
\begin{align}\label{compare-2}
v_2^k(y)=V_k(y)+M_kH_{1,y}(Q_1^k)\bigg (\frac{\pi}4|\frac{p_k}{\epsilon_k M_k}|^2M_k-2\pi\frac{\Delta \mathfrak{h}_2^k(0)}{\mathfrak{h}_2^k(0)}\frac{\delta_k^2}{M_k}\bigg )\\
+16\pi \epsilon_k(\partial_1H_{1,y}(Q_1^k)c_1+\partial_2 H_{1,y}(Q_1^k)c_2)\bigg )+o(M_k\epsilon_k). \nonumber
\end{align}
From (\ref{compare-1}) and (\ref{compare-2}) we have
\begin{align}\label{compare-3}
    &\frac{\bar V_k-V_k}{M_k}
    =(H_{1,y}(Q_1^k)-H_{2,y}(e_1))\bigg (\frac{\pi}{4}|\frac{p_k}{M_k\epsilon_k}|^2 M_k -2\pi\frac{\Delta \mathfrak{h}_2^k(0)}{\mathfrak{h}_2^k(0)}\frac{\delta_k^2}{M_k}\bigg ) \\
    +&16\pi\epsilon_k((\partial_1H_{1,y}(Q_1^k)+\partial_1H_{2,y}(e_1))c_1
    +(\partial_2 H_{1,y}(Q_1^k)+\partial_2 H_{2,y}(e_1))c_2)+o(\epsilon_k)
    \nonumber
\end{align}
$$H_{1,y}(Q_1^k)-H_{2,y}(e_1)=G_k(y,Q_1^k)-G_k(y,e_1). $$
Now we evaluate the left hand side when $y$ is away from $e_1$ and $e^{i\pi}$.
\begin{align*}
&\bar V_k(y)-V_k(y)\\
=&2\log \frac{1+\frac{e^{\mu_k}}{D_k}|y^2-1|^2}{1+\frac{e^{\mu_k}}{D_k}|y^2-1-p_k|^2}\\
=&2\log \frac{e^{-\mu_k}D_k+|y^2-1|^2}{e^{-\mu_k}D_k+|y^2-1-p_k|^2}.
\end{align*}
To evaluate the second term we have
$$|y^2-1-p_k|^2=|y^2-1|^2(1-2Re(\frac{y^2-1}{|y^2-1|^2}\bar p_k)+
\frac{|p_k|^2}{|y^2-1|^2}). $$
In order for the comparison in (\ref{compare-3}), obviously $M_k=O(\epsilon_k)$. This fact means $|p_k|=O(\epsilon_k^2)$. The the left hand side of (\ref{compare-3}) is 
$$\frac{\bar V_k(y)-V_k(y)}{M_k}=4Re(\frac{y^2-1}{|y^2-1|^2}\frac{\bar p_k}{M_k\epsilon_k})\epsilon_k+O(\epsilon_k^3).$$
Then 
$$(H_{1,y}(Q_1^k)-H_{2,y}(e_1))\tilde M_k=-\frac{1}{2\pi}\tilde M_k\log \frac{|y-1|}{|y+1|}+o(\epsilon_k)$$
where 
$$\tilde M_k=\frac{\pi}{4}(\frac{|p_k|}{M_k\epsilon_k})^2M_k-2\pi\frac{\Delta\mathfrak{h}_2^k(0)}{\mathfrak{h}_2^k(0)}\frac{\delta_k^2}{M_k}. $$
$$c_1^2+c_2^2=|\frac{p_k}{4M_k\epsilon_k}|^2. $$
The second term on the right hand side is
$$2\epsilon_k\bigg ((\frac{y_1+1}{|y+1|^2}+\frac{y_1-1}{|y-1|^2}+1)\cos\theta_1+(\frac{y_2}{|y+1|^2}+\frac{y_2}{|y-1|^2})\sin\theta_1\bigg ).$$
If $M_k\ge C \epsilon_k$, the $O(\delta_k^2/M_k)$ is an error term.  By choosing $|y|$ large we see that the logarithmic term is close to $1$, while other terms tend to zero, which violates the equality. If  $M_k=o(\epsilon_k)$, (\ref{compare-3}) is also violated if $|y|$ is large because the terms of $y$ tend to zero in different rates.
Lemma \ref{small-other} is established. $\Box$

Proposition \ref{key-w8-8} is an immediate consequence of Lemma \ref{small-other}.  $\Box$.

\medskip

Now we finish the proof of (\ref{vanish-first-tau}).

Let $\hat w_k=w_k/\tilde \delta_k$. (Recall that $\tilde \delta_k=\delta_k |\nabla \mathfrak{h}_2^k(0)|+\delta_k^2\mu_2^k$). If $|\nabla \mathfrak{h}_2^k(0)|/(\delta_k\mu_2^k)\to \infty$, we see that in this case 
$\tilde \delta_k\sim \delta_k\mu_2^k|\nabla \mathfrak{h}_2^k(0)|$. The equation of $\hat w_k$ is
\begin{equation}\label{hat-w}
\Delta \hat w_k+|y|^{2}e^{\xi_k}\hat w_k=a_k\cdot (e_1-y)|y|^{2}e^{V_k}+b_ke^{V_k}|y-e_1|^{2}|y|^{2}e^{V_k},
\end{equation}
in $\Omega_k$, where $a_k=\delta_k\nabla \mathfrak{h}_2^k(0)/\tilde \delta_k$, $b_k=o(1/\mu_2^k)$. 
By Proposition \ref{key-w8-8}, $|\hat w_k(y)|\le C$. Before we carry out the remaining part of the proof we observe that $\hat w_k$ converges to a harmonic function in $\mathbb R^2$ minus finite singular points. Since $\hat w_k$ is bounded, all these singularities are removable. Thus $\hat w_k$ converges to a constant. Based on the information around $e_1$, we shall prove that this constant is $0$. However, looking at the right hand side the equation,
$$ (e_1-y)|y|^{2}e^{V_k}\rightharpoonup 8\pi  (e_1-e^{i\pi})\delta_{e^{i\pi}}. $$
we will get a contradiction by comparing the Pohozaev identities of $v_2^k$ and $V_k$, respectively.

Now we use the notation $W_k$ again and use Proposition \ref{key-w8-8} to rewrite the equation for $W_k$.
Let
$$W_k(z)=\hat w_k(e_1+\epsilon_k z), \quad |z|< \delta_0 \epsilon_k^{-1} $$
for $\delta_0>0$ small. Then from Proposition \ref{key-w8-8} we have
\begin{equation}\label{h-exp}
\mathfrak{h}_2^k(\delta_ky)=\mathfrak{h}_2^k(\delta_k e_1)+\delta_k \nabla \mathfrak{h}_2^k(\delta_k e_1)(y-e_1)+O(\delta_k^2)|y-e_1|^2,
\end{equation}
\begin{equation}\label{y-1}
|y|^{2}=|e_1+\epsilon_k z|^{2}=1+O(\epsilon_k)|z|,
\end{equation}
\begin{equation}\label{v-radial}
V_k(e_1+\epsilon_k z)+2\log \epsilon_k=U_k(z)+O(\epsilon_k)|z|+O(\epsilon_k^2)(\log (1+|z|))^2
\end{equation}
and
\begin{equation}\label{xi-radial}
\xi_k(e_1+\epsilon_k z)+2\log \epsilon_k=U_k(z)+O(\epsilon_k)(1+|z|).
\end{equation}
Using (\ref{h-exp}),(\ref{y-1}),(\ref{v-radial}) and (\ref{xi-radial}) in (\ref{hat-w}) we write the equation of $W_k$ as
\begin{equation}\label{W-rough}
\Delta W_k+\mathfrak{h}_2^k(\delta_k e_1)e^{U_k(z)}W_k=-\epsilon_k a_k\cdot ze^{U_k(z)}+E_w, \quad 0<|z|<\delta_0 \epsilon_k^{-1}
\end{equation}
where
\begin{equation}\label{ew-rough-2}
E_w(z)=O(\epsilon_k)(1+|z|)^{-3}, \quad |z|<\delta_0 \epsilon_k^{-1}.
\end{equation}

Since $\hat w_k$ obviously converges to a global harmonic function with removable singularity, we have $\hat w_k\to \bar c$ for some $\bar c\in \mathbb R$. Then we claim that
\begin{lem}\label{bar-c-0} $\bar c=0$.
\end{lem}

\noindent{\bf Proof of Lemma \ref{bar-c-0}:}

 If $\bar c\neq 0$, we use $W_k(z)=\bar c+o(1)$ on $B(0,\delta_0 \epsilon_k^{-1})\setminus B(0, \frac 12\delta_0 \epsilon_k^{-1})$ and consider the projection of $W_k$ on $1$:
 $$g_0(r)=\frac 1{2\pi}\int_{0}^{2\pi}W_k(re^{i\theta})d\theta. $$
If we use $F_0$ to denote the projection to $1$ of the right hand side we have, using the rough estimate of $E_w$ in (\ref{ew-rough-2})
$$g_0''(r)+\frac 1r g_0'(r)+\mathfrak{h}_2^k(\delta_ke_1)e^{U_k(r)}g_0(r)=F_0,\quad 0<r<\delta_0 \epsilon_k^{-1} $$
where
$$F_0(r)=O(\epsilon_k)(1+|z|)^{-3}. $$
In addition we also have
$$\lim_{k\to \infty} g_0(\delta_0 \epsilon_k^{-1})=\bar c+o(1). $$
For simplicity we omit $k$ in some notations. By the same argument as in Lemma \ref{w-around-e1},  we have
$$g_0(r)=O(\epsilon_k)\log (2+r),\quad 0<r<\delta_0 \epsilon_k^{-1}. $$
Thus $\bar c=0$.
Lemma \ref{bar-c-0} is established. $\Box$

\medskip

Based on Lemma \ref{bar-c-0} and standard Harnack inequality for elliptic equations we have
\begin{equation}\label{small-til-w}
\tilde w_k(x)=o(1),\,\,\nabla \tilde w_k(x)=o(1),\,\, x\in B_3\setminus B(e^{i\pi},\delta_0).
\end{equation}
Equation (\ref{small-til-w}) is equivalent to $w_k=o(\tilde \delta_k)$ and $\nabla w_k=o(\tilde \delta_k)$ in the same region.

\medskip

In the next step we consider the difference between two Pohozaev identities.
We consider the Pohozaev identity around $Q_1^k$. Let $\Omega_{1,k}=B(Q_1^k,\tau)$ for small $\tau>0$. For $v_k$ we have
(\ref{pi-vk}). Correspondingly the Pohozaev identity for $V_k$ is

\begin{align*}
\int_{\Omega_{1,k}}\partial_{\xi}(|y|^{2}\mathfrak{h}_2^k(\delta_ke_1))e^{V_k}-\int_{\partial \Omega_{1,k}}e^{V_k}|y|^{2}\mathfrak{h}_2^k(\delta_k e_1)(\xi\cdot \nu)\\
=\int_{\partial \Omega_{1,k}}(\partial_{\nu}V_k\partial_{\xi}V_k-\frac 12|\nabla V_k|^2(\xi\cdot \nu))dS. \nonumber
\end{align*}
Using $w_k=v_2^k-V_k$ and $|w_k(y)|\le C\tilde \delta_k$ we have
\begin{align*}
&\int_{\partial \Omega_{1,k}}(\partial_{\nu}v_2^k\partial_{\xi}v_2^k-\frac 12|\nabla v_2^k|^2(\xi\cdot \nu))dS\\
=&\int_{\partial \Omega_{1,k}}(\partial_{\nu}V_k\partial_{\xi}V_k-\frac 12|\nabla V_k|^2(\xi\cdot \nu))dS\\
&+\int_{\partial \Omega_{1,k}}(\partial_{\nu}V_k\partial_{\xi}w_k+\partial_{\nu}w_k\partial_{\xi}V_k-(\nabla V_k\cdot \nabla w_k)(\xi\cdot \nu))dS+o(\tilde \delta_k).
\end{align*}
If we just use crude estimate: $\nabla w_k=o(\tilde \delta_k)$, we have
\begin{align*}&\int_{\partial \Omega_{1,k}}(\partial_{\nu}v_2^k\partial_{\xi}v_2^k-\frac 12|\nabla v_2^k|^2(\xi\cdot \nu))dS\\
-&\int_{\partial \Omega_{1,k}}(\partial_{\nu}V_k\partial_{\xi}V_k-\frac 12|\nabla V_k|^2(\xi\cdot \nu))dS
=o(\tilde \delta_k).
\end{align*}
The difference on the second terms is minor: If we use the expansion of $v_2^k=V_k+w_k$ and that of $\mathfrak{h}_2^k(\delta_ky)$ around $e_1$, it is easy to obtain
$$\int_{\partial \Omega_{1,k}}e^{v_k}|y|^{2}\mathfrak{h}_2^k(\delta_ky)(\xi\cdot \nu)-\int_{\partial \Omega_{1,k}}e^{V_k}|y|^{2}\mathfrak{h}_2^k(\delta_ke_1)(\xi\cdot \nu)=o(\tilde \delta_k). $$
To evaluate the first term,  we use
\begin{align}\label{imp-1}
&\partial_{\xi}(|y|^{2}\mathfrak{h}_2^k(\delta_ky))e^{v_2^k}\\
=&\partial_{\xi}(|y|^{2}\mathfrak{h}_2^k(\delta_ke_1)+|y|^{2}\delta_k\nabla \mathfrak{h}_2^k(\delta_ke_1)(y-e_1)+O(\delta_k^2))e^{V_k}(1+w_k+O(\delta_k^2\mu_k))\nonumber\\
=&\partial_{\xi}(|y|^{2})\mathfrak{h}_2^k(\delta_k e_1)e^{V_k}+\delta_k\partial_{\xi}(|y|^{2}\nabla \mathfrak{h}_2^k(\delta_ke_1)(y-e_1))e^{V_k}\nonumber\\
&+\partial_{\xi}(|y|^{2}\mathfrak{h}_2^k(\delta_ke_1))e^{V_k}w_k+O(\delta_k^2\mu_k)e^{V_k}.\nonumber
\end{align}

For the third term on the right hand side of (\ref{imp-1}) we use the equation for $w_k$:
$$\Delta w_k+\mathfrak{h}_2^k(\delta_ke_1)e^{V_k}|y|^{2}w_k=-\delta_k\nabla \mathfrak{h}_2^k(\delta_ke_1)\cdot (y-e_1)|y|^{2}e^{V_k}+O(\delta_k^2)e^{V_k}|y|^{2}.$$
From integration by parts we have
\begin{align}\label{extra-1}
&\int_{\Omega_{1,k}}\partial_{\xi}(|y|^{2})\mathfrak{h}_2^k(\delta_ke_1)e^{V_k}w_k\nonumber\\
=&2\int_{\Omega_{1,k}}y_{\xi}\mathfrak{h}_2^k(\delta_ke_1)e^{V_k}w_k\nonumber\\
=&2\int_{\Omega_{1,k}}\frac{y_{\xi}}{|y|^2}(-\Delta w_k-\delta_k \nabla\mathfrak{h}_2^k(\delta_ke_1)(y-e_1)|y|^{2}e^{V_k}+O(\delta_k^2)e^{V_k}|y|^{2})\nonumber\\
=&-2\delta_k\int_{\Omega_{1,k}}\frac{y_{\xi}}{|y|^{2}}\nabla \mathfrak{h}_2^k(\delta_ke_1)(y-e_1)|y|^{2}e^{V_k}\nonumber\\
&+2\int_{\partial \Omega_{1,k}}(\partial_{\nu}(\frac{y_{\xi}}{|y|^2})w_k-\partial_{\nu}w_k\frac{y_{\xi}}{|y|^2})+o(\tilde \delta_k)\nonumber\\
=&\nabla\mathfrak{h}_2^k(\delta_ke_1)\bigg (-16\delta_k\pi(e^{i\pi}\cdot \xi)(e^{i\pi}-e_1)+O(\mu_2^k\epsilon_k^2)\bigg )+o(\tilde \delta_k),
\end{align}
where we have used $\nabla w_k, w_k=o(\tilde \delta_k)$ on $\partial \Omega_{s,k}$.
For the second term on the right hand side of (\ref{imp-1}), we have
\begin{align}\label{imp-2}
&\int_{\Omega_{1,k}}\delta_k\partial_{\xi}(|y|^{2}\nabla \mathfrak{h}_2^k(\delta_ke_1)(y-e_1))e^{V_k}\\
=&2\delta_k\int_{\Omega_{1,k}}y_{\xi}\nabla \mathfrak{h}_2^k(\delta_ke_1)(y-e_1)e^{V_k}+
\delta_k\int_{\Omega_{1,k}}|y|^{2}\partial_{\xi}\mathfrak{h}_2^k(\delta_ke_1)e^{V_k} \nonumber \\
=&\nabla \mathfrak{h}_2^k(\delta_ke_1)\big (16N\pi\delta_k(e^{i\pi}\cdot \xi)(e^{i\pi}-e_1)+O(\mu_k\epsilon_k^2)\big )\nonumber\\
&+
\delta_k\partial_{\xi}\mathfrak{h}_2^k(\delta_ke_1)(8\pi+O(\mu_2^k\epsilon_k^2))+o(\tilde \delta_k).\nonumber
\end{align}
Using (\ref{extra-1}) and (\ref{imp-2}) in the difference between (\ref{pi-vk}) and (\ref{pi-Vk}), we have
$$\delta_k \partial_{\xi}\mathfrak{h}_2^k(\delta_ke_1)(1+O(\mu_2^k\epsilon_k^2))=o(\tilde \delta_k). $$
(\ref{vanish-first-tau}) is established.

\subsection{Laplacian Vanishing Property}

First we consider the case that $\delta_k\le C\mu_2^k\epsilon_k$. In this case $\epsilon_k^{-1}\delta_k^2\le C\epsilon_k^{\epsilon}$ for some $\epsilon\in (0,1)$. The whole argument of Proposition \ref{key-w8-8} can be
employed to prove
\begin{equation}\label{2nd-w}
|w_k(y)|\le C\delta_k^2(\mu_2^k)^{\frac{7}{4}} 
\end{equation}
In order to employ the same strategy of proof, one needs to have two things: first $\epsilon_k^{-1}\delta_k^2=O(\epsilon_k^{\epsilon})$. This is clear from the definition of $\delta_k$. Second, in the proof of Lemma \ref{small-other} and in (\ref{subtle-1}) we need 
$$\delta_k^{3}/M_k=o(\epsilon_k), $$
where $M_k>\delta_k^2(\mu_2^k)^{7/4}$. Since $\delta_k\le C\mu_2^k\epsilon_k$, the required inequality holds.  The proof of Proposition \ref{key-w8-8} follows.

\medskip
 
The precise upper bound of $w_k$ in (\ref{2nd-w}) leads to the vanishing rate of the Laplacian estimate : If we use 
$$W_k(z)=w_k(e^{i\pi}+\epsilon_kz)/(\delta_k^2(\mu_2^k)^{\frac 74}), \quad |z|<\tau \epsilon_k^{-1}$$ where $e_l\neq e_1$. 
We shall show that the projection of $W_k$ over $1$ is not bounded when $|z|\sim \epsilon_k^{-1}$, which gives the desired contradiction. 

We write the equation of $w_k$ as
$$\Delta w_k+|y|^{2}e^{\xi_k}w_k
=(\mathfrak{h}_2^k(\delta_k e_1)-\mathfrak{h}_2^k(\delta_k y))|y|^{2}e^{V_k}. 
$$
Then 
\begin{align*}
&\Delta W_k(z)+e^{U_k}W_k(z)\\
=&a_0e^{U_k}+a_1ze^{U_k}+\frac 1{2\mu_k^{7/4}}\Delta \mathfrak{h}_2^k(0)|z|^2e^{U_k}+\frac{1}{\mu_k^{7/4}}R_2(\theta)|z|^2e^{U_k}
+O(\epsilon_k^{\epsilon}(1+|z|)^{-3}).
\end{align*}
where 
$$a_0=(\mathfrak{h}_2^k(\delta_ke_1)-\mathfrak{h}_2^k(\delta_ke_l))/(\delta_k^2(\mu_2^k)^{7/4}),$$
$$a_1=-\nabla\mathfrak{h}_2^k(\delta_ke_l)/(\delta_k(\mu_2^k)^{7/4}),$$
$R_2$ is the collection of spherical harmonic functions of degree $2$. Note that there is no appearance of $\epsilon_k$ or $\epsilon_k^2$ in the equation for $W_k$. This is a key point in the proof. Around $e_1$, the terms are small, while around $-1$, the terms are far greater. This discrepancy leads to contradiction.

Let $g_k(r)$ be the projection of $W_k$ on $1$, by the same ODE analysis as before, we see that $g_k$ satisfies
$$g_k''+\frac 1rg_k'(r)+e^{U_k}g_k=E_k $$
where
$$E_k(r)=O(\epsilon_k^{\epsilon})(1+r)^{-3}+\frac 1{2(\mu_2^k)^{7/4}}\Delta (\log \mathfrak{h}_2^k)(0)r^2e^{U_k}. $$

Using the same argument as in Lemma \ref{w-around-e1}, we have
 $$g_k(r)\sim \Delta (\log \mathfrak{h}_2^k)(0)
(\log r)^2 \mu_k^{-7/4},\quad r>10. $$

Clearly if $\Delta (\log \mathfrak{h}_2^k(0)\neq 0$ we obtain a violation of the bound of $w_k$ for $r\sim \epsilon_k^{-1}$. The laplacian vanishing theorems is proved under the assumption 
\begin{equation}\label{assump-small}
\epsilon_k^{-1}|Q_1^k-e^{i\pi}|\le \epsilon_k^{\epsilon}.
\end{equation}
We need this assumption because the $\xi_k$ function that comes from the equation of $w_k$ needs to tend to $U$ after scaling. From (3.13) in \cite{wei-zhang-adv}, 
$|Q_1^k-e^{i\pi}|=O(\delta_k^2)+O(\mu_2^ke^{-\mu_2^k})$.
If $\delta_k^2\epsilon_k^{-1}\ge C$, the argument in the proof of (\ref{vanish-first-tau}) cannot be used because either $\xi_k$ does not tend $U$ or $c_0=0$ cannot be proved.

\medskip

\noindent{\bf Proof of the Laplace Vanishing property for $\delta_k\ge \mu_2^k\epsilon_k$. }

In this case we  write the equation of $w_k$ as
$$\Delta w_k+|y|^2\mathfrak{h}_2^k(\delta_ky)e^{v_k}-|y|^2\mathfrak{h}_2^k(\delta_ke_1)e^{V_k}=0. $$

From $0=\nabla w_k(e_1)$ we have
\begin{equation}\label{location-mid-1}
0=\int_{\Omega_k}\nabla_1G_k(e_1,\eta)|\eta |^2(\mathfrak{h}_2^k(\delta_k \eta)e^{v_k}-\mathfrak{h}_2^k(\delta_ke_1)e^{V_k})d\eta+O(\delta_k^3)
\end{equation}
Note that $v_k$ is close to another global solution $\bar V_k$ which matches with a local maximum of $v_k$ at $Q_2^k$. 
Evaluating the right hand side of (\ref{location-mid-1}) we have
$$\nabla_1 G_k(e_1,Q_2^k)-\nabla_1 G_k(e_1,e^{i\pi})=O(\epsilon_k^2\mu_k)+O(\delta_k^3). $$
This expression gives
$$Q_2^k-e^{i\pi}=O(\delta_k^3)+O(\mu_k\epsilon_k^2). $$
This estimate will lead to a better estimate of $w_k$ outside the two bubbling disks. 
From the Green's representation for $w_k$ we now obtain 
$$w_k(y)=\int_{\Omega_k}(G_k(y,\eta)-G_k(e_1,\eta))|\eta |^2(\mathfrak{h}_2^k(\delta_k\eta)e^{v_2^k(y)}-\mathfrak{h}_2^k(\delta_ke_1)e^{V_k})d\eta+O(\delta_k^2)$$
where the last term $O(\delta_k^2)$ comes from the oscillation of $w_k$ on $\partial \Omega_k$. Then we have
\begin{align*}
    w_k(y)
    &=-\frac 1{2\pi}\int_{\Omega_k}\log \frac{|y-\eta |}{ |e_1-\eta |}|\eta |^2(\mathfrak{h}_2^k(\delta_k \eta)e^{v_2^k}-\mathfrak{h}_2^k(\delta_k e_1)e^{V_k})+O(\delta_k^2)\\
    &=-4\log \frac{|y-Q_2^k|}{|e_1-Q_2^k|}+4\log \frac{|y-e^{i\pi}|}{2}+O(\delta_k^2\mu_k).
\end{align*}
By $|Q_2^k-e^{i\pi}|=O(\delta_k^2)$ we see that $w_k(y)=O(\delta_k^2)$ on $|y-e^{i\pi}|=\tau$.

The standard point-wise estimate for singular equation ( see \cite{zhangcmp,gluck} ) gives 
\begin{align*}
&v_k(Q_2^k+\epsilon_kz)+2\log \epsilon_k\\
=&\log \frac{e^{\mu_2^k}}{(1+\frac{e^{\mu_2^k}}{8\mathfrak{h}_2^k(\delta_kQ_k)}|z|^2)^2}+\phi_1^k+C\delta_k^2\Delta (\log \mathfrak{h}_2^k)(0)(\log (1+|z|))^2, \quad |z|\sim \epsilon_k^{-1}. 
\end{align*}

\begin{align*}
&V_k(e^{i\pi}+\epsilon_kz)+2\log \epsilon_k \\
=&\log \frac{e^{\mu_2^k}}{(1+\frac{e^{\mu_2^k}}{8\mathfrak{h}_2^k(\delta_ke_1)}|z|^2)^2}+\phi_2^k+O(\epsilon_k^2(\log \epsilon_k)^2),\quad |z|\sim \epsilon_k^{-1}. 
\end{align*}

Thus 
$$w_k(Q_2^k+\epsilon_kz)=O(\epsilon_k^2(\log \epsilon_k)^2)+\phi_1^k-\phi_2^k+C\Delta (\log \mathfrak{h}_2^k)(0)\delta_k^2(\log (1+|z|)^2), $$
for $ |z|\sim \epsilon_k^{-1}. $
Taking the average around the origin, the spherical averages of the two harmonic functions are zero and $O(\delta_k^2)$ respectively, since they take zero at the origin and a point  at most $O(\delta_k^2)$ from the origin. So the spherical average of $w_k$ is comparable to 
$$\Delta (\log \mathfrak{h}_2^k)(0)\delta_k^2 (\log \epsilon_k)^2$$ for $|z|\sim \epsilon_k^{-1}$. Thus we know $\Delta (\log \mathfrak{h}_2^k)(0)=o(1)$ because $w_k=O(\delta_k^2\mu_k)$ in this region, 
The Laplacian vanishing theorem is established for all the cases. Theorem \ref{local-laplace} is established. $\Box$

\medskip

Obviously the proof of Theorem \ref{local-laplace} can be used to rule out $(4,2)$ and $(2,4)$ if the spherical Harnack inequality is violated. 

\section{Proof of Theorem \ref{main-thm}}

  Let $u^k$ be a sequence of blowup solutions and let $p$ be a blowup point. In a local coordinate, if around $p$ there is a $(2,4)$ type blowup profile with a violation of the spherical Harnack inequality around the local maximum of $u_1^k$, then we see from the
  derivation of (\ref{main-local}) that $u^k=(u_1^k,u_2^k)$ can be written as a locally defined blowup solutions: $\mathfrak{u}^k=(\mathfrak{u}_1^k,\mathfrak{u}_2^k)$ of 
  \begin{align*}
     \Delta \mathfrak{u}_1^k+2\mathfrak{h}_1^ke^{\mathfrak{u}_1^k}-\mathfrak{h}_2^ke^{\mathfrak{u}_2^k}=0, \\
     \Delta \mathfrak{u}_2^k-\mathfrak{h}_1^ke^{\mathfrak{u}_1^k}+2\mathfrak{h}_2^ke^{\mathfrak{u}_2^k}=0, 
  \end{align*}
  where 
  $$\mathfrak{h}_i^k(x)=\rho_i^kh_i^k(x)e^{\phi_i+f_i}, \quad i=1,2 $$
  where $\phi_i$ and $f_i$ are defined in (\ref{local-phi}) and (\ref{local-f}), respectively. 
 Since in this case we have $\lim_{k\to \infty} \Delta (\log \mathfrak{h}_2^k)(0)=0$ and 
$$\Delta (\log \mathfrak{h}_2^k)(0)=\Delta \log h_2^k(p_k)-2K(p_k)+4m\pi $$
for some $m\in \mathbb Z$ we obtain a contradiction to the curvature assumption (\ref{curvature-a}). The assumption that $\mu_1^k>\mu_2^k+5\log \delta_k^{-1}$ is guaranteed by the existence of a fully bubbling blowup point. Around that blowup point, $u_1^k-u_2^k=O(1)$ away from bubbling disks. Then the crude estimate of $u_1^k$ and $u_2^k$ around the three disks in the formation $(2,4)$ gives that $\mu_1^k=\mu_2^k+(6+o(1))\log \delta_k^{-1}+O(1)$. Thus the two types $(2,4)$ and $(4,2)$ are ruled out if their formations violate the spherical Harnack inequality.

In order to prove the main theorem we also need to rule out the blowup type $(4,4)$ if the spherical Harnack inequality is violated. In the formation of $(4,4)$. There are four situations where the spherical Harnack inequality is violated. Because of similarly we only describe two of them. In the first case there are four bubbling disks all tending to one point. Among the four disks, two bubbling disks of $u_2^k$ are in the same line of one bubbling disk of $u_1^k$ and the distance from each local maximum of $u_2^k$ to the local maximum of $u_1^k$ is $\delta_k$. On the other hand there is another bubbling disk of $u_1^k$ whose center (which is also a local maximum of $u_1^k$) is $\bar \delta_k$ away from the previous local maximum of $u_1^k$. In this case $\bar \delta_k/\delta_k\to \infty$ even though $\bar \delta_k\to 0$. 

The second case is that there are two bubbling disks, in the first bubbling disk, $u_1^k$ and $u_2^k$ both satisfy spherical Harnack inequality and the height of $u_1^k$ minus the height of $u_2^k$ tends to infinity. If the equation for $u_2^k$ is scaled according to the height of $u_2^k$, the new function $v_2^k$ after scaling tends to a global solution of
$$\Delta v_2+|y|^2e^{v_2}=4\pi \delta_p,\quad \mbox{in}\quad \mathbb R^2. $$
Beside this bubbling disk, there is another disk of $u_1^k$ that is of distance $\delta_k$ away from the center of the first bubbling disk.

Now we rule out the first case: Let $p_0^k$ be the center of the disk of $u_1^k$ that is placed in the middle of two $u_2^k$ disks. By setting
$$\tilde v_i^k(y)=u_i^k(p_0^k+\delta_ky)+2\log \delta_k,\quad i=1,2. $$
For this function we set $\Omega_k=B(0,\tau \delta_k^{-1})$. 
Note that there is another blowup point $p_1^k\to p$, but $|p_1^k-p_k|/\delta_k\to \infty$. We use $\bar \delta_k:=|p_1^k-p_0^k|$. So after scaling, the $p_1^k$ becomes
$L_k=\bar \delta_k/\delta_k$ away from the origin. We first give a rough description of the heights of bubbling disks. Let $\tilde v_1^k(p_0)=\mu_1^k$, let $\mu_2^k$ be the height of $\tilde v_2^k$ in one of the two disks. Let $\bar \mu_1^k$ be the height of $\tilde v_1^k$ at the other bubbling disk. By calculating the spherical average around $0$ and $Q_k$ (which is the image of $p_1^k$ after scaling, $|Q_k|\sim L_k$), we have
$$-\mu_1^k=-\bar\mu_1^k-4\log L_k+o(1)L_k+O(1). $$
Because of the fully bubbling disk, $u_1^k-u_2^k=O(1)$ on $\partial B(p_0^k,\tau)$, this leads to 
$$-\mu_2^k-(6+o(1))\log L_k=-\bar \mu_1^k-(4+o(1))\log L_k+O(1). $$
Thus we have 
$$\bar \mu_1^k=\mu_2^k+(2+o(1))\log L_k$$
$$\mu_1^k=\mu_2^k+(6+o(1))\log L_k. $$
In the setting of previous sections we use the following function $f_1^k$ to remove $\tilde v_1^k$ from the equation for $\tilde v_2^k$:
$$f_1^k(y)=-\int_{\Omega_k}G_k(y,\eta)\mathfrak{h}_1^k(\delta_k \cdot)e^{\tilde v_1^k}d\eta+\int_{\Omega_k}G_k(0,\eta)\mathfrak{h}_1^k(\delta_k \cdot)e^{\tilde v_1^k}d\eta. $$
So $f_1^k$ satisfies 
$$\Delta f_1^k=\mathfrak{h}_1^k(\delta_k y),\quad f_1^k(0)=0. $$
For $|y|\sim 1$ we have 
$$f_1^k(y)=2\log |y|+2\log |y-Q_k|+o(\delta_k^2). $$
Applying the gradient estimate as before, we first have 
$$\delta_k\nabla \mathfrak{h}_2^k(0)-2\frac{Q_k}{|Q_k|^2} =O(\mu_2^ke^{-\mu_2^k})+O(\delta_k^2)$$
Since $\frac{1}{|Q_k|}=\delta_k/\bar \delta_k$, the equality above cannot hold if $\mu_2^ke^{-\mu_2^k}=O(\delta_k)$. Thus we consider $\mu_2^ke^{-\mu_2^k}>>\delta_k$, which implies $\mu_2^k=O(\log \delta_k^{-1})$. In this case the second gradient estimate gives
$$\delta_k\nabla \mathfrak{h}_2^k(0)-\frac{2 Q_k}{|Q_k|^2}=O(\delta_k^2\mu_2^k). $$
But this cannot hold because the $\frac{Q_k}{|Q_k|^2}$ majorizes all other terms. 
We have ruled out the case of $(4,4)$ if the spherical Harnack inequality does not hold in the formation of $(2,4)$.

\medskip

Next we rule out the case that $(4,4)$ consists of one bubbling disk of $(2,4)$ and $(2,0)$ but the $(2,4)$ satisfies a spherical Harnack around some point. 
In other words $u_1^k$ stays on top of $u_2^k$. We use $\bar \lambda_1^k$ to represent the height of $u_1^k$ in the bubbling disk that has local maximums of $u_2^k$, we let $\lambda_2^k$ represent the height of $u_2^k$ and set $\lambda_1^k$ to represent $u_1^k$ in the other disk. The distance between these two bubbling disks is $\delta_k\to 0$. 

Let $p_k$ be a local maximum of $u_2^k$ and $q_k$ be the local maximum of $u_1^k$, which is also the center of a bubbling disk outside the bubbling disk of $u_2^k$. By the context just mentioned, we set $2\delta_k=|p_k-q_k|$. If we set
$$v_i^k(y)=u_i^k(p_k+\epsilon_{2,k}y)+2\log \epsilon_{2,k}, \quad i=1,2. $$
where $\epsilon_{2,k}=e^{-\lambda_{2,k}/2}$, then
$v_2^k$ converges, along a sub-sequence to 
$$\Delta v_2+2e^{v_2}=4\pi\delta_0,\quad \mbox{in}\quad \mathbb R^2. $$
Note that the limit function $v_2$ also satisfies $\int_{\mathbb R^2}e^{v_2}<\infty$ because of the restriction from $u_2^k$. 
By the classification theorem of Prajapat-Tarantello \cite{prajapat}
$$v_2(y)=-6\log |y|+O(1),\quad |y|>1. $$
Then for $|x-p_k|\sim \delta_k$ we have
$$u_2^k(y)=-6\log \delta_k+4\log \epsilon_{2,k}+O(1)=-6\log \delta_k-2\lambda_2^k+O(1). $$

Since one assumption of Theorem \ref{main-thm} is that there is a fully bubbling point. The estimate of fully bubbling sequence, as well as Harnack inequality imply that for $|x-p_k|\sim \tau$, $u_1^k-u_2^k=O(1)$. By Green's representation formulas one can also obtain the finite difference for $$|u_1^k(x)-u_2^k(x)|\le C, 
\quad |x-p_k|=5\delta_k. $$

If we employ the point-wise estimate for single Liouville equation (see Li \cite{li-cmp}) we have, for $|x-q_k|\sim \delta_k$, 
$$ u_1^k(x)=-\lambda_1^k-4\log \delta_k+O(1),\quad |x-q_k|=\delta_k/2. $$
On the other hand, for $v_2^k$ we have
$$v_2^k(y)=-6\log \frac{\delta_k}{\epsilon_{2,k}}+O(1),\quad |y|=\frac 12 \delta_k/\epsilon_{2,k}, $$
which is equivalent to 
$$u_2^k(x)=-2\lambda_2^k-6\log \delta_k+O(1),\quad |x-p_k|=\frac 12\delta_k. $$
Thus we have
$$\lambda_1^k=2\lambda_2^k+2\log \delta_k+O(1). $$
Note that the appearance of bubbles implies
$$\lambda_2^k+2\log \delta_k\to \infty $$

Next we derive the relation between $\bar \lambda_1^k$ and $\lambda_1^k$. 

The equation of $\frac 13 v_1^k+\frac 23v_2^k$ is
$$ -\Delta (\frac 13 v_1^k+\frac 23 v_2^k)=\mathfrak{h}_2^k(p_k+\epsilon_{2,k}y)e^{v_2^k(y)},
\quad |y|\le \frac 12\delta_k/\epsilon_{2,k}. $$

Because of the fast decay of $e^{v_2^k}$ it is easy to use a bounded function to remove it from the right hand side of the equation. Then we have
\begin{equation}\label{add-eq-1}
(\frac 13 v_1^k+\frac 23v_2^k)|_{\partial B_1}=(\frac 13 v_1^k+\frac 23v_2^k)|_{\partial B_{L_k}}+O(1) 
\end{equation}
where $L_k=\frac 12\delta_k\epsilon_{2,k}^{-1}$. 
The value of $v_1^k$ on $\partial B_1$ is $-(\bar \lambda_1^k-\lambda_2^k)+O(1)$ based on standard estimates for single Liouville equation ( \cite{li-cmp}). It is also obvious that 
$v_2^k=O(1)$ on $\partial B_1$. On the other hand $v_1^k=-\lambda_1^k-4\log \delta_k-\lambda_2^k+O(1)$ on $\partial B_{L_k}$ based on the analysis of single equation. In fact it is also the value for $v_2^k$ because of the influence of the fully bubbling profile. If we compute the value of $v_2^k$ from its limit, it is
$$v_2^k(y)=-6\log (\delta_k/\epsilon_{2,k})+O(1),\quad |y|=L_k. $$
Using these information in (\ref{add-eq-1}) we have
$$\bar \lambda_1^k=\lambda_2^k+9(\lambda_2^k+2\log \delta_k)+O(1). $$

Now we consider
$$\hat v_1^k(y)=u_1^k(\epsilon_{2,k}y+p_k)+2\log \epsilon_{2,k}. $$
The Pohozaev identity gives
$$\epsilon_{2,k}\nabla\log h_1^k(p_k)+\nabla \phi_k(0)=O(\epsilon_3^k\mu_3^k) $$
where $\phi_k$ is the harmonic function that eliminates the oscillation of $\hat v_1^k$ around $B_1$. 
$$\epsilon_{3,k}=e^{-\frac{\bar \lambda_1^k-\lambda_2^k}2},\quad \mu_3^k=\bar \lambda_1^k-\lambda_2^k=9(\lambda_2^k+2\log \delta_k)+O(1). $$
Recall that
\begin{equation}\label{lambda-12-d}
\bar \lambda_1^k-\lambda_2^k=9(\lambda_2^k+2\log \delta_k)+O(1). \end{equation}
The harmonic function $\phi_k$ is determined by the singular source, which is now located at $\frac{\delta_k}{\epsilon_{2,k}}$ away from the origin. The derivative of this function is comparable to $\epsilon_{2,k}/\delta_k$.
So the contradiction boils down to whether
$$\frac{\epsilon_{2,k}}{\delta_k}>>\epsilon_{3,k}^2\mu_{3,k}. $$
Using (\ref{lambda-12-d})
we see that we need 
$$e^{4(\lambda_2^k+2\log \delta_k)}>>9(\lambda_2^k+2\log\delta_k). $$
Since this quantity goes to infinity, it holds.
Theorem \ref{main-thm} is established. $\Box$

\section{Appendix A}
In this section we provide a proof of (\ref{late-1}). 
Here we briefly explain the roles of each term. $\phi_1$ corresponds to the radial solution in the kernel of the linearized operator of the global equation. In other words, $\phi_1^k/\bar M_k$ should tend to zero. $\phi_2^k/\bar M_k$ is the combination of the two other functions in the kernel. $\phi_4$ is the second order term which will play a leading role later. $\phi_3^k$ comes from the difference of $\mathfrak{h}_2^k$ at $Q_1^k$ and $e_1$. Before we derive (\ref{late-1}) we point out that 
$$\frac{|p_k|}{\epsilon_k}\le C\epsilon_k^{-1}\mu_2^ke^{-\mu_2^k}+C\epsilon_k^{-1}\delta_k^2\le C\mu_2^k\epsilon_k. $$

The derivation of (\ref{late-1}) is as follows: First by the expression of $e_1$ in (\ref{ql-exp}) we have 
$$y^{2}=1+2\epsilon_kz+\epsilon_k^2z^2, $$
where $y=e_1+\epsilon_k z$. Then
$$|y^2-e_1|^2=4\epsilon_k^2|z+\frac{\epsilon_k}2z^2|^2 $$

\begin{align*}
&|y^{2}-e_1-p_k|^2\\
=&4\epsilon_k^2\bigg (|z|^2-2Re(\frac{\epsilon_k}2|z|^2\bar z)-2Re(z\frac{\bar p_k}{2\epsilon_k})+|\frac{p_k}{2\epsilon_k}|^2
+O(\epsilon_k^2|z|^4)+O(\mu_2^k\epsilon_k^2|z|^2)\bigg ).
\end{align*}
Next by the definition of $D_k$ in (\ref{def-D}) 
$$\frac{D_k-\bar D_k}{\bar D_k}=\delta_k \nabla (\log \mathfrak{h}_2^k)(\delta_k e_1)\cdot (Q_1^k-1)+O(\delta_k^2).$$
\begin{align}\label{tem-5}
\frac{e^{\mu_2^k-\bar \mu_2^k}}{\bar D_k}
&=\frac{1}{D_k}(1+\frac{D_k-\bar D_k}{\bar D_k}+\bar \mu_2^k-\mu_2^k+O(\bar \mu_2^k-\mu_2^k)^2+O(\delta_k^2)). \\
&=\frac 1{D_k}(1+\delta_k\nabla\log \mathfrak{h}_2^k(\delta_k e_1)\cdot (Q_1^k-e_1)+\bar \mu_2^k-\mu_2^k+O(\bar \mu_2^k-\mu_2^k)^2+O(\delta_k^2)). \nonumber
\end{align}
Then the expression of $A$ is ( for simplicity we omit $k$ in some notations)
\begin{align*}
    A=\bigg (\frac{e^{\bar \mu_2-\mu_2}}{\bar D_k}(4\big (|z|^2-2Re(z\frac{\overline{p_k}}{2\epsilon_k})+\frac{|p_k|^2}{4\epsilon_k^2}+O(\epsilon_k|z|^3)\big )\\
    -\frac{4}{D_k}(|z|^2+O(\epsilon_k)|z|^3)\bigg )/B
\end{align*}
After using (\ref{tem-5}) we have 
\begin{align}\label{exp-A-2}
A=\bigg (\frac{1}{D_k}(\delta_k\nabla(\log \mathfrak{h}_2^k(\delta_ke_1)(Q_1-1)+\bar\mu_2^k-\mu_2^k+O(\bar \mu_2^k-\mu_2^k)^2)4|z|^2\\
-4Re(z\frac{\bar p_k}{\epsilon_k})\frac 1{D_k}+
\frac{|p_k|^2}{\epsilon_k^2D_k}+O(\epsilon_k)|z|^3+O(\delta_k^2)|z|^2\bigg )/B. \nonumber
\end{align}

\begin{equation}\label{exp-a-squ}
A^2=
(\frac{16}{D_k^2}(Re(z\frac{-\bar p_k}{\epsilon_k})^2)/B^2
+\mbox{other terms}.
\end{equation}
The numerator of $A^2$ has the following leading term:
$$\frac{8}{D_k^2}\bigg (|z|^2(\frac{|p_k|}{\epsilon_k})^2\big (1+\cos(2\theta-2\theta_{1}^k)\big ) \bigg )$$
where $z=|z|e^{i\theta}$, $p_k=|p_k|e^{i\theta_1^k}$.
Using these expressions we can obtain (\ref{late-1}) by direct computation. 

\section{Appendix B}
In this section we derive (\ref{late-2-r}). Here we first point out that we shall use (\ref{late-1}).
The terms of $\phi_1$ and $\phi_3$ will lead to $o(\epsilon_k)$, the integration involving $\phi_2$ cancels with the second term of (\ref{late-2-r}). The computation of $\phi_2$ is based on this equation:
$$\int_{\mathbb R^2}\frac{\frac{\mathfrak{h}_2^k(\delta_k Q_1^k)}4\sigma_k\nabla \mathfrak{h}_2^k(\delta_k Q_1^k)(Q_1^k-e_1)|z|^2}{(1+\frac{\mathfrak{h}_2^k(\delta_k Q_1^k)}8|z|^2)^3}dz
=8\pi \sigma_k\nabla (\log \mathfrak{h}_2^k)(\delta_k Q_1^k)(Q_1^k-e_1), $$
and by (\ref{delta-small-1})
\begin{equation}\label{subtle-1}
\nabla \log \mathfrak{h}_2^k(\delta_kQ_1^k)-\nabla \log \mathfrak{h}_2^k(\delta_k e_1)=O(\delta_k)=o(\epsilon_k). 
\end{equation}

The integration involving $\phi_4$ provides the leading term. More detailed information is the following:
First for a global solution
$$V_{\mu,p}=\log \frac{e^{\mu}}{(1+\frac{e^{\mu}}{\lambda}|z^{N+1}-p|^2)^2}$$ of
$$\Delta V_{\mu,p}+\frac{8(N+1)^2}{\lambda}|z|^{2N}e^{V_{\mu,p}}=0,\quad \mbox{in }\quad \mathbb R^2, $$
by differentiation with respect to $\mu$ we have
$$\Delta(\partial_{\mu}V_{\mu,p})+\frac{8(N+1)^2}{\lambda}|z|^{2N}e^{V_{\mu,p}}\partial_{\mu}V_{\mu,p}=0,\quad \mbox{in}\quad \mathbb R^2. $$
By the expression of $V_{\mu,p}$ we see that 
$$\partial_r\bigg (\partial_{\mu}V_{\mu,p}\bigg )(x)=O(|x|^{-2N-3}).$$ 
Thus we have
\begin{align}\label{inte-eq-1}
&\int_{\mathbb R^2}\partial_{\mu}V_{\mu,p}|z|^{2N}e^{V_{\mu,p}}\\
=&\int_{\mathbb R^2}\frac{(1-\frac{e^{\mu}}{\lambda}|z^{N+1}-P|^2)|z|^{2N}}{(1+\frac{e^{\mu}}{\lambda}|z^{N+1}-P|^2)^3}dz\nonumber \\
=&-\frac{\lambda}{8(N+1)^2}\int_{\mathbb R^2}\Delta (\partial_{\nu}V_{\mu,p})=0. \nonumber
\end{align}

From $V_{\mu,p}$  we also have
$$\int_{\mathbb R^2}\partial_{P}V_{\mu,p}|y|^{2N}e^{V_{\mu,p}}=\int_{\mathbb R^2}\partial_{\bar P}V_{\mu,p}|y|^{2N}e^{V_{\mu,p}}=0, $$
which gives
\begin{equation}\label{inte-eq-2}
\int_{\mathbb R^2}\frac{\frac{e^{\mu}}{\lambda}(\bar z^{N+1}-\bar P)|z|^{2N}}{(1+\frac{e^{\mu}}{\lambda}|z^{N+1}-P|^2)^3}=\int_{\mathbb R^2}\frac{\frac{e^{\mu}}{\lambda}( z^{N+1}- P)|z|^{2N}}{(1+\frac{e^{\mu}}{\lambda}|z^{N+1}-P|^2)^3}=0.
\end{equation}

Now we need more precise expressions of $\phi_1$, $\phi_3$ and $B$:
\begin{align*}
&\phi_1=(\mu_2^k-\bar \mu_2^k)(1-\frac{4}{D_k}|z+\frac{1}2\epsilon_k z^2|^2)/B, \\
&\phi_3=\frac{8}{D_k B}Re((z+\frac{1}2\epsilon_kz^2))(\frac{-\bar p_k}{\epsilon_k}))
\\
&B=1+\frac{4}{D_k}|z+z^2\epsilon_k|^2,
\end{align*}
From here we use scaling and cancellation to have
$$\int_{B(0,\tau\epsilon_k^{-1})}\frac{\phi_1}{M_k}B^{-2}=o(\epsilon_k),$$
which is based on (\ref{inte-eq-1}) and
$$\int_{B(0,\tau\epsilon_k^{-1})}\frac{\phi_3}{M_k}B^{-2}=o(\epsilon_k)$$
which is based on (\ref{inte-eq-2}). 
Thus (\ref{late-2-r}) holds.

\end{document}